\documentclass[11pt]{amsart}
\usepackage[dvipsnames]{xcolor}
\usepackage{tikz}
\usepackage{scalerel}
\usepackage{marvosym} 
\usepackage[margin=2.4cm]{geometry}
\usepackage{float}
\allowdisplaybreaks
\newif\ifverbose
\verbosetrue 

\newcommand{\drawchessboard}[2]{%
\begin{tikzpicture}[scale=0.23]
  \foreach \x in {-14,...,13} {
    \foreach \y in {-14,...,13} {
      \pgfmathtruncatemacro{\isblack}{mod(\x+\y,2)==0 ? 1 : 0}
      \ifnum\isblack=1
        \fill[gray!25] (\x,\y) rectangle ++(1,1);
      \fi
    }
  }
  #2
  \draw[thick, red] (-5,-5) rectangle (5,5);
   \draw[thick, gray] (-14,-14) rectangle (14,14);
\end{tikzpicture}%
}

\newcommand{\chessboardtwoQueen}[2]{%
\begin{tikzpicture}[scale=0.26]
  \foreach \x in {-5,...,4} {
    \foreach \y in {-5,...,4} {
      \pgfmathtruncatemacro{\isblack}{mod(\x+\y,2)==0 ? 1 : 0}
      \ifnum\isblack=1
        \fill[gray!25] (\x,\y) rectangle ++(1,1);
      \fi
    }
  }
  #2
  \draw[thick, red] (-2,-2) rectangle (2,2);
   \draw[thick, gray] (-5,-5) rectangle (5,5);
\end{tikzpicture}%
}

\newcommand{\chessboardthreeQueen}[2]{%
\begin{tikzpicture}[scale=0.26]
  \foreach \x in {-6,...,5} {
    \foreach \y in {-6,...,5} {
      \pgfmathtruncatemacro{\isblack}{mod(\x+\y,2)==0 ? 1 : 0}
      \ifnum\isblack=1
        \fill[gray!25] (\x,\y) rectangle ++(1,1);
      \fi
    }
  }
  #2
  \draw[thick, red] (-2,-2) rectangle (2,2);
   \draw[thick, gray] (-6,-6) rectangle (6,6);
\end{tikzpicture}%
}

\newcommand{\chessboardthreeQueenodd}[2]{%
\begin{tikzpicture}[scale=0.24]
  \foreach \x in {-6,...,6} {
    \foreach \y in {-6,...,6} {
      \pgfmathtruncatemacro{\isblack}{mod(\x+\y,2)==0 ? 1 : 0}
      \ifnum\isblack=1
        \fill[gray!25] (\x,\y) rectangle ++(1,1);
      \fi
    }
  }
  #2
  \draw[thick, red] (-2,-2) rectangle (3,3);
   \draw[thick, gray] (-6,-6) rectangle (7,7);
\end{tikzpicture}%
}

\newcommand{\chessboardfourQueeneven}[2]{%
\begin{tikzpicture}[scale=0.22]
  \foreach \x in {-8,...,7} {
    \foreach \y in {-8,...,7} {
      \pgfmathtruncatemacro{\isblack}{mod(\x+\y,2)==0 ? 1 : 0}
      \ifnum\isblack=1
        \fill[gray!25] (\x,\y) rectangle ++(1,1);
      \fi
    }
  }
  #2
  \draw[thick, red] (-3,-3) rectangle (3,3);
   \draw[thick, gray] (-8,-8) rectangle (8,8);
\end{tikzpicture}%
}

\newcommand{\chessboardfourQueenodd}[2]{%
\begin{tikzpicture}[scale=0.22]
  \foreach \x in {-8,...,8} {
    \foreach \y in {-8,...,8} {
      \pgfmathtruncatemacro{\isblack}{mod(\x+\y,2)==0 ? 1 : 0}
      \ifnum\isblack=1
        \fill[gray!25] (\x,\y) rectangle ++(1,1);
      \fi
    }
  }
  #2
  \draw[thick, red] (-2,-2) rectangle (3,3);
   \draw[thick, gray] (-8,-8) rectangle (9,9);
\end{tikzpicture}%
}

\newcommand{\chessboardfiveQueen}[2]{%
\begin{tikzpicture}[scale=0.22]
  \foreach \x in {-9,...,8} {
    \foreach \y in {-9,...,8} {
      \pgfmathtruncatemacro{\isblack}{mod(\x+\y,2)==0 ? 1 : 0}
      \ifnum\isblack=1
        \fill[gray!25] (\x,\y) rectangle ++(1,1);
      \fi
    }
  }
  #2
  \draw[thick, red] (-3,-3) rectangle (3,3);
   \draw[thick, gray] (-9,-9) rectangle (9,9);
\end{tikzpicture}%
}

\newcommand{\chessboardfiveQueenodd}[2]{%
\begin{tikzpicture}[scale=0.22]
  \foreach \x in {-9,...,9} {
    \foreach \y in {-9,...,9} {
      \pgfmathtruncatemacro{\isblack}{mod(\x+\y,2)==0 ? 1 : 0}
      \ifnum\isblack=1
        \fill[gray!25] (\x,\y) rectangle ++(1,1);
      \fi
    }
  }
  #2
  \draw[thick, red] (-2,-2) rectangle (3,3);
   \draw[thick, gray] (-9,-9) rectangle (10,10);
\end{tikzpicture}%
}

\newcommand{\chessboardsixQueen}[2]{%
\begin{tikzpicture}[scale=0.22]
  \foreach \x in {-10,...,10} {
    \foreach \y in {-10,...,10} {
      \pgfmathtruncatemacro{\isblack}{mod(\x+\y,2)==0 ? 1 : 0}
      \ifnum\isblack=1
        \fill[gray!25] (\x,\y) rectangle ++(1,1);
      \fi
    }
  }
  #2
  \draw[thick, red] (-4,-4) rectangle (5,5);
   \draw[thick, gray] (-10,-10) rectangle (11,11);
\end{tikzpicture}%
}

\newcommand{\chessboardsixQueeneven}[2]{%
\begin{tikzpicture}[scale=0.22]
  \foreach \x in {-11,...,10} {
    \foreach \y in {-11,...,10} {
      \pgfmathtruncatemacro{\isblack}{mod(\x+\y,2)==0 ? 1 : 0}
      \ifnum\isblack=1
        \fill[gray!25] (\x,\y) rectangle ++(1,1);
      \fi
    }
  }
  #2
  \draw[thick, red] (-4,-4) rectangle (4,4);
   \draw[thick, gray] (-11,-11) rectangle (11,11);
\end{tikzpicture}%
}

\newcommand{\chessboardsevenQueen}[2]{%
\begin{tikzpicture}[scale=#1]
  \foreach \x in {-12,...,12} {
    \foreach \y in {-12,...,12} {
      \pgfmathtruncatemacro{\isblack}{mod(\x+\y,2)==0 ? 1 : 0}
      \ifnum\isblack=1
        \fill[gray!25] (\x,\y) rectangle ++(1,1);
      \fi
    }
  }
  #2
  \draw[thick, red] (-3,-3) rectangle (4,4);
   \draw[thick, gray] (-12,-12) rectangle (13,13);
\end{tikzpicture}%
}

\newcommand{\chessboardsevenQueeneven}[2]{%
\begin{tikzpicture}[scale=#1]
  \foreach \x in {-13,...,12} {
    \foreach \y in {-13,...,12} {
      \pgfmathtruncatemacro{\isblack}{mod(\x+\y,2)==0 ? 1 : 0}
      \ifnum\isblack=1
        \fill[gray!25] (\x,\y) rectangle ++(1,1);
      \fi
    }
  }
  #2
  \draw[thick, red] (-4,-4) rectangle (4,4);
   \draw[thick, gray] (-13,-13) rectangle (13,13);
\end{tikzpicture}%
}

\newcommand{\chessboardnineQueenodd}[2]{%
\begin{tikzpicture}[scale=0.17]
  \foreach \x in {-15,...,15} {
    \foreach \y in {-15,...,15} {
      \pgfmathtruncatemacro{\isblack}{mod(\x+\y,2)==0 ? 1 : 0}
      \ifnum\isblack=1
        \fill[gray!25] (\x,\y) rectangle ++(1,1);
      \fi
    }
  }
  #2
  \draw[thick, red] (-4,-4) rectangle (5,5);
   \draw[thick, gray] (-15,-15) rectangle (16,16);
\end{tikzpicture}%
}

\usepackage{amsmath,amssymb}
\usepackage{enumerate}
\usepackage{skak}
\usepackage{tikz}
\usetikzlibrary{arrows.meta,automata,quotes}
\usetikzlibrary{backgrounds}
\usetikzlibrary{matrix}
\usepackage{pgfplots}

\usepackage{microtype}
\usepackage[colorlinks]{hyperref}
\usepackage[T1]{fontenc}
\usepackage[utf8]{inputenc}
\usepackage{mathrsfs}  
\usepackage{caption}
\theoremstyle{plain}
\newtheorem{thm}{Theorem}
\newtheorem{cor}[thm]{Corollary}
\newtheorem{prop}[thm]{Proposition}
\newtheorem{lem}[thm]{Lemma}

\newtheorem{thm*}{Theorem}
 
\theoremstyle{definition}
\newtheorem{defi} [thm] {Definition}

\newtheorem{remark}[thm]{Remark}

\newtheorem{obs}[thm]{Observation}

\newcommand{\Z}{\mathbb{Z}}
\newcommand{\N}{\ensuremath{\mathbb{N}}}

\newcommand{\n}{\ensuremath{\ensuremath{\left\{\left\lfloor\frac{2-n}{2}\right\rfloor,\ldots, \left\lfloor\frac{n}{2}\right\rfloor\right\}}}}
\newcommand{\B}{\mathcal B}
\newcommand{\C}{\mathcal C}
\newcommand{\R}{\mathcal R}
\renewcommand{\O}{\mathcal O}
\renewcommand{\S}{\mathcal S}
\renewcommand{\ge}{\geqslant}
\renewcommand{\le}{\leqslant}

\renewcommand{\leq}{\leqslant}
\newcommand{\inl}{\operatorname{inloss}}
\newcommand{\cl}{\operatorname{cenloss}}
\newcommand{\loss}{\operatorname{loss}}
\newcommand{\cover}{\operatorname{cover}}

\title{Thresholds of Queen covers}
\author{Harman Agrawal, Sahana Jahagirdar, Prem Kant, Urban Larsson}
\address{Dept. IEOR, IIT Bombay}
\email{harman.agrawal@gmail.com, sahanaj@iitb.ac.in, premkant072@gmail.com, larsson@iitb.ac.in}
\date{\today}

\begin{document}
\begin{abstract}
    We study optimal configurations of Queens on a square chessboard, defined as those covering the maximum number of squares. For a fixed number of Queens, $q$, we prove the existence of two thresholds in board size: a non-attacking threshold beyond which all optimal configurations are pairwise non-attacking, and a stabilizing threshold beyond which the set of optimal configurations becomes constant. Related studies on Queen domination, such as Tarnai and Gáspár (2007), focus on minimizing the number of Queens needed for full board coverage. Our approach, by contrast, fixes the number of Queens and analyzes optimal cover via a certain loss-function due to {\em internal loss} and {\em decentralization}. We demonstrate how the internal loss can be decomposed in terms of defined concepts, {\em balance} and {\em overlap concentration}. Moreover, by using our results, for sufficiently large board sizes, we find all optimal Queen configurations for all $2\le q\le 9$. And, whenever possible, we relate those solutions in terms of the classical problem of placing $q$ non-attacking Queens on a $q\times q$ board. For example, in case $q=8$, out of the twelve classical fundamental solutions, only three apply here as centralized patterns on large boards. On the other hand, the single classical fundamental solution for $q=6$ is never cover optimal on large boards, even if centralized, but another pattern that fits inside a $q\times (q+1)$ board applies.
\end{abstract}
\maketitle

\section{Introduction}\label{sec:intro} 

How many Queens of Chess are required to cover every square of a board of size $n$ by $n$?  Each Queen {\em covers} her own position, as well as all orthogonal and diagonal squares that she can reach in one move. This question has been researched in the literature (e.g. \cite{AV2017, BM2002, TG2007, W2002}), and nice surveys appear in \cite{C1990, JW2004}.\footnote{The classical $8$-Queens problem, and more generally the ``$N$-Queens problem'', is often credited to Max Bezzel \cite{MB1850}, but it is reported that Carl Friedrich Gauss also solved it, finding 72 solutions. This account appears in \cite{EL1883}, based on a letter from Gauss’s son Eugene. See also \cite{W} for a brief survey of this problem. For a recent related study on infinite board sizes, see \cite{DS2019}. } If $n=3$, then a centralized Queen covers the full board. If $n=4$, the problem is already more interesting. Only two Queens are required, unless we impose that no Queen can attack any other Queen. Perhaps counter-intuitively, in that case, it is only possible to cover the full board with three Queens. In a sense, the crux lies in `the size of the board' and the requirement that the whole board must be covered. We illustrate these motivating examples in Figure~\ref{fig:4by4}.

\begin{figure} [ht!]
\begin{center}
\begin{tikzpicture} [scale = 0.36]
\foreach \x in {-1,...,1} {
    \foreach \y in {-1,...,1} {
        \pgfmathparse{mod(\x+\y,2)==0 ? "gray!15" : "white"}
        \edef\cellcolor{\pgfmathresult}
        \fill[fill=\cellcolor] (\x, \y) rectangle ++(1,1);
    }
}
\foreach \x/\y in {0/0} 
    \draw (\x+.5, \y+.5) node {${\tiny \symqueen}$};


\def\off{0.3}

\draw[red!60, thick] (0.5+\off,0.5+\off) -- (1.5,1.5);   
\draw[red!60, thick] (0.5-\off,0.5-\off) -- (-0.5,-0.5); 
\draw[red!60, thick] (0.5-\off,0.5+\off) -- (-0.5,1.5);  
\draw[red!60, thick] (0.5+\off,0.5-\off) -- (1.5,-0.5);  

\draw[red!60, thick] (0.5-\off,0.5) -- (-0.5,0.5); 
\draw[red!60, thick] (0.5+\off,0.5) -- (1.5,0.5);  
\draw[red!60, thick] (0.5,0.5-\off) -- (0.5,-0.5); 
\draw[red!60, thick] (0.5,0.5+\off) -- (0.5,1.5);  

\draw[step=1cm,gray, thin] (-1, -1) grid (2,2); 
\end{tikzpicture}
\hspace{5 mm}
\begin{tikzpicture} [scale = 0.36]
\foreach \x in {-2,...,1} {
    \foreach \y in {-2,...,1} {
        \pgfmathparse{mod(\x+\y,2)==0 ? "gray!15" : "white"}
        \edef\cellcolor{\pgfmathresult}
        \fill[fill=\cellcolor] (\x, \y) rectangle ++(1,1);
    }
}
\foreach \x/\y in {-1/0,0/0} 
    \draw (\x+.5, \y+.5) node {${\tiny \symqueen}$};
\draw[step=1cm,gray, thin] (-2, -2) grid (2,2); 
\end{tikzpicture}\hspace{5 mm}
\begin{tikzpicture} [scale = 0.36]
\foreach \x in {-2,...,1} {
    \foreach \y in {-2,...,1} {
        \pgfmathparse{mod(\x+\y,2)==0 ? "gray!15" : "white"}
        \edef\cellcolor{\pgfmathresult}
        \fill[fill=\cellcolor] (\x, \y) rectangle ++(1,1);
    }
}
\foreach \x/\y in {0/-2,-1/1,1/0} 
    \draw (\x+.5, \y+.5) node {${\tiny \symqueen}$};
\draw[step=1cm,gray, thin] (-2, -2) grid (2,2); 
\end{tikzpicture}
\end{center}
\caption{A centralized Queen covers the whole $3\times 3$ board. The $4\times 4$ board can be covered by two Queens, but three are required if they must be non-attacking.}
\label{fig:4by4}
\end{figure}

In this paper we relax the problem to, for a given finite number of Queens, finding {\em optimal} solutions, i.e. those maximizing coverage. We prove that: 
\begin{enumerate}
    \item for a given number of Queens, $q$, there is a threshold (that may depend on $q$) on the size of the square board such that, beyond that threshold every optimal configuration of Queen contains only non-attacking Queens. 
    \item for a given number of Queens, $q$, there is a threshold on the size of the board such that, beyond that threshold every optimal configuration of Queens  belongs to a constant finite set of optimal Queen configurations. 
\end{enumerate}
The first (second) threshold will be called the non-attacking (stabilizing) threshold; see Theorems~\ref{thm:main1} and \ref{thm:main2}, respectively. Let us illustrate by depicting some boards with four Queen configurations. Similar to the case with two Queens on a $4\times 4$ board, on a $9\times 9$ board there exists an optimal configuration with a pair of attacking Queens. However, on larger boards, optimality implies that no Queen attacks any other Queen. In fact, the set of optimal representatives for a $10\times 10$ board has size two modulo symmetry; see Figure~\ref{fig:4Qs}. But, by exhaustive search, already on an $11\times 11$ board, another configuration outperforms these two candidates. 
In Figure~\ref{fig:4Qs}, we show some configurations on boards of side lengths 10 and 12. We use  dark purple shading to illustrate winning board configurations. It turns out that the winning pattern for the $12\times 12$  board remains optimal for all larger boards although the second stabilizing threshold occurs at a slightly larger board size. We discuss the exact value of the stabilizing threshold in Section~\ref{sec:app}.

We will demonstrate that optimality of cover corresponds to minimizing a defined {\em loss}, which, for large boards become constant, and this is the approach of this paper.  
\begin{figure} [ht!]
\definecolor{boarddark}{HTML}{8c198c}
\definecolor{boardlight}{HTML}{D8BFE6}
\definecolor{blocked}{HTML}{F2F2F2}
\definecolor{queen}{HTML}{2C3E50}
\definecolor{queencrown}{HTML}{C9A227}

\begin{center}
\begin{tikzpicture}[scale = 0.36]

\foreach \x in {-5,...,4} {
    \foreach \y in {-5,...,4} {
        \pgfmathsetmacro{\covered}{
            (\x==-4) || (\y==-2) || (abs(\x+4)==abs(\y+2)) || 
            (\x==-2) || (\y==3)  || (abs(\x+2)==abs(\y-3))  || 
            (\x==1)  || (\y==-4) || (abs(\x-1)==abs(\y+4))  || 
            (\x==3)  || (\y==1)  || (abs(\x-3)==abs(\y-1))     
        }
        \ifdim \covered pt>0pt
            \fill[boarddark] (\x,\y) rectangle ++(1,1);
        \fi
    }
}
\foreach \x/\y in {-4/-2, -2/3, 1/-4, 3/1}
 \node at (\x+0.5,\y+0.5) {\textcolor{yellow}{${\small \symqueen}$}};

\draw[step=1cm,black,thin] (-5,-5) grid (5,5);
\end{tikzpicture}
\hspace{5 mm}
\begin{tikzpicture}[scale = 0.36]
\foreach \x in {-5,...,4} {
    \foreach \y in {-5,...,4} {
        \pgfmathsetmacro{\covered}{
            (\x==-3) || (\y==-2) || (abs(\x+3)==abs(\y+2)) || 
            (\x==-2) || (\y==2)  || (abs(\x+2)==abs(\y-2))  || 
            (\x==1)  || (\y==-3) || (abs(\x-1)==abs(\y+3))  || 
            (\x==2)  || (\y==1)  || (abs(\x-2)==abs(\y-1))     
        }
        \ifdim \covered pt>0pt
            \fill[boarddark] (\x,\y) rectangle ++(1,1);
        \fi
    }
}

\foreach \x/\y in {-3/-2, -2/2, 1/-3, 2/1}
   \node at (\x+0.5,\y+0.5) {\textcolor{yellow}{${\small \symqueen}$}};

\draw[step=1cm,black,thin] (-5,-5) grid (5,5);
\end{tikzpicture}
\hspace{5 mm}
\begin{tikzpicture}[scale = 0.36]
\foreach \x in {-5,...,4} {
    \foreach \y in {-5,...,4} {
        \pgfmathsetmacro{\covered}{
            (\x==0) || (\y==-2) || (abs(\x+0)==abs(\y+2)) || 
            (\x==-2) || (\y==-1)  || (abs(\x+2)==abs(\y+1))  || 
            (\x==1)  || (\y==0) || (abs(\x-1)==abs(\y+0))  || 
            (\x==-1)  || (\y==1)  || (abs(\x+1)==abs(\y-1))     
        }
        \ifdim \covered pt>0pt
            \fill[boardlight] (\x,\y) rectangle ++(1,1);
        \fi
    }
}

\foreach \x/\y in {0/-2,-1/1,-2/-1,1/0} 
    \node at (\x+0.5,\y+0.5) {${\small \symqueen}$};

\draw[step=1cm,black,thin] (-5,-5) grid (5,5);
\end{tikzpicture}
\vspace{5 mm}

\begin{tikzpicture}[scale = 0.36]
\foreach \x in {-6,...,5} {
    \foreach \y in {-6,...,5} {
        \pgfmathsetmacro{\covered}{
            (\x==-4) || (\y==-2) || (abs(\x+4)==abs(\y+2)) || 
            (\x==-2) || (\y==3)  || (abs(\x+2)==abs(\y-3))  || 
            (\x==1)  || (\y==-4) || (abs(\x-1)==abs(\y+4))  || 
            (\x==3)  || (\y==1)  || (abs(\x-3)==abs(\y-1))     
        }
        \ifdim \covered pt>0pt
            \fill[boardlight] (\x,\y) rectangle ++(1,1);
        \fi
    }
}

\foreach \x/\y in {-4/-2, -2/3, 1/-4, 3/1}
    \node at (\x+0.5,\y+0.5) {${\small \symqueen}$};

\draw[step=1cm,black,thin] (-6,-6) grid (6,6);
\end{tikzpicture}
\hspace{5 mm}
\begin{tikzpicture}[scale = 0.36]
\foreach \x in {-6,...,5} {
    \foreach \y in {-6,...,5} {
        \pgfmathsetmacro{\covered}{
            (\x==-3) || (\y==-2) || (abs(\x+3)==abs(\y+2)) || 
            (\x==-2) || (\y==2)  || (abs(\x+2)==abs(\y-2))  || 
            (\x==1)  || (\y==-3) || (abs(\x-1)==abs(\y+3))  || 
            (\x==2)  || (\y==1)  || (abs(\x-2)==abs(\y-1))     
        }
        \ifdim \covered pt>0pt
            \fill[boardlight] (\x,\y) rectangle ++(1,1);
        \fi
    }
}


\foreach \x/\y in {-3/-2, -2/2, 1/-3, 2/1}
    \node at (\x+0.5,\y+0.5) {${\small \symqueen}$};

\draw[step=1cm,black,thin] (-6,-6) grid (6,6);
\end{tikzpicture}
\hspace{5 mm}
\begin{tikzpicture}[scale = 0.36]
\foreach \x in {-6,...,5} {
    \foreach \y in {-6,...,5} {
        \pgfmathsetmacro{\covered}{
            (\x==0) || (\y==-2) || (abs(\x+0)==abs(\y+2)) || 
            (\x==-2) || (\y==-1)  || (abs(\x+2)==abs(\y+1))  || 
            (\x==1)  || (\y==0) || (abs(\x-1)==abs(\y+0))  || 
            (\x==-1)  || (\y==1)  || (abs(\x+1)==abs(\y-1))     
        }
        \ifdim \covered pt>0pt
            \fill[boarddark] (\x,\y) rectangle ++(1,1);
        \fi
    }
}
\foreach \x/\y in {0/-2,-1/1,-2/-1,1/0} 
\node at (\x+0.5,\y+0.5) {\textcolor{yellow}{${\small \symqueen}$}};

\draw[step=1cm,black,thin] (-6,-6) grid (6,6);
\end{tikzpicture}
\end{center}
\caption{Three configurations of four non-attacking Queens for sidelengths $n=10$ and $n=12$ respectively; in each picture white color indicates non-covered squares. Yellow Queens on a dark board indicate optimal configurations. The two upper leftmost pictures are the only optimal configurations for a $10\times 10$ board (modulo mirror symmetry).
For $12\times 12$, however, the configuration at the lower right is optimal and outperforms those to the left.}
\label{fig:4Qs}
\end{figure}

In Section~\ref{sec:loss}, we define the setting and relate cover to the concept of {\em loss}, via a combination of defined {\em internal} and {\em centralized loss}. In Section~\ref{sec:nonattack} we prove the existence of the non-attacking threshold and in Section~\ref{sec:stablize} we establish existence of the stabilizing threshold. 
In Section~\ref{sec:balovercon}, the more intriguing internal loss is dissected via the notions of {\em balance} and {\em overlap concentration}. In Section~\ref{sec:app}, we list results on the stabilizing threshold for $2\le q\le 9$. 
In Section~\ref{sec:history} we view our results from a historical perspective, and in Section~\ref{sec:fut}, we mention some open problems and future research directions. Finally, in Section~\ref{sec:thoughts}, we share some thoughts about further applications of this study.

\section{Computation of Queen cover via loss}\label{sec:loss} 
Let us review the setup. It is convenient to centralize the board with $(0,0)$ in the middle, if the side length is odd, and otherwise $(0,0)$ is the lower left square of the four central ones. Thus the $n\times n$ board, $n\in \N=\{0,1,\ldots \}$, is
$$\B_n=\n\times\n.$$ 
For example, $\B_0=\varnothing$, $\B_1=\{(0,0)\}$, $\B_2=\{(0,0),(0,1),(1,0),(1,1)\}$ and $$\B_3=\{(-1,-1),(-1,0),(-1,1),(0,-1),(0,0),(0,1),(1,-1),(1,0),(1,1)\};$$ see also Figure~\ref{fig:visualization}. We skip the axes for simplicity and write $0$ instead to identify coordinates; as in  Figure~\ref{fig:loss} onwards.

\begin{figure}[ht!]
\centering
\begin{tikzpicture}[scale=0.36]

 \draw[gray, thin] (-0.5,-0.5) rectangle (0.5,0.5);       

\draw[<->] (-3,0) -- (4,0) node[right] {$~$};
\draw[<->] (0,-3) -- (0,4) node[above] {$~$};

\foreach \x in {-2,...,3}{
  \draw (\x,0.1) -- (\x,-0.1);
}

\foreach \y in {-2,...,3}{
  \draw (0.1,\y) -- (-0.1,\y);
}

\foreach \x in {0,0}
    \foreach \y in {0,0}
        \filldraw[black] (\x,\y) circle (2pt);

\end{tikzpicture}
\hspace{5mm}
\begin{tikzpicture}[scale=0.36]

\draw[gray, thin] (-0.5,-0.5) rectangle (1.5,1.5);       
\foreach \x in {0.5,0.5}                                
    \draw[gray, thin] (\x,-0.5) -- (\x,1.5);
\foreach \y in {0.5,0.5}                                
    \draw[gray, thin] (-0.5,\y) -- (1.5,\y);

\draw[<->] (-3,0) -- (4,0) node[right] {$~$};
\draw[<->] (0,-3) -- (0,4) node[above] {$~$};

\foreach \x in {-2,...,3}{
  \draw (\x,0.1) -- (\x,-0.1);
}

\foreach \y in {-2,...,3}{
  \draw (0.1,\y) -- (-0.1,\y);
}

\foreach \x in {0,1}
    \foreach \y in {0,1}
        \filldraw[black] (\x,\y) circle (2pt);

\end{tikzpicture}
\hspace{5mm}
\begin{tikzpicture}[scale=0.36]


\draw[gray,thin] (-1.5,-1.5) rectangle (1.5,1.5);

\foreach \x in {-0.5,0.5}
  \draw[gray,thin] (\x,-1.5) -- (\x,1.5);
\foreach \y in {-0.5,0.5}
  \draw[gray,thin] (-1.5,\y) -- (1.5,\y);

\draw[<->] (-3,0) -- (4,0) node[right] {$~$};
\draw[<->] (0,-3) -- (0,4) node[above] {$~$};

\foreach \x in {-2,...,3}{
  \draw (\x,0.1) -- (\x,-0.1);
}

\foreach \y in {-2,...,3}{
  \draw (0.1,\y) -- (-0.1,\y);
}
\foreach \x in {-1,0,1}
  \foreach \y in {-1,0,1}
    \filldraw[black] (\x,\y) circle (2pt);  

\end{tikzpicture}
\hspace{5mm}
\begin{tikzpicture}[scale=0.36]


\draw[gray,thin] (-1.5,-1.5) rectangle (2.5,2.5);

\foreach \x in {-0.5,0.5,1.5}
  \draw[gray,thin] (\x,-1.5) -- (\x,2.5);
\foreach \y in {-0.5,0.5,1.5}
  \draw[gray,thin] (-1.5,\y) -- (2.5,\y);

\draw[<->] (-3,0) -- (4,0) node[right] {$~$};
\draw[<->] (0,-3) -- (0,4) node[above] {$~$};

\foreach \x in {-2,...,3}{
  \draw (\x,0.1) -- (\x,-0.1);
}

\foreach \y in {-2,...,3}{
  \draw (0.1,\y) -- (-0.1,\y);
}

\foreach \x in {-1,0,1,2}
  \foreach \y in {-1,0,1,2}
    \filldraw[black] (\x,\y) circle (2pt);  

\end{tikzpicture}
\caption{A visualization of the boards $\B_1$, $\B_2$, $\B_3$, and $\B_4$.}\label{fig:visualization}
\end{figure}

A {\em configuration} of Queens is a finite set $\C\subset \Z\times \Z$. A configuration is a $q$-configuration if it contains exactly $q$ Queens. An attacking space,  
$$A:2^{\Z\times \Z}\rightarrow 2^{\Z\times \Z},$$ 
is  a function on the configurations that outputs the total set of attacked squares on a doubly infinite grid. That is, for any configuration $\C$, 
$$A(\C)=\{(x+i,y), (x,y+i), (x+i,y+i), (x+i,y-i)\mid (x,y)\in \C, i\in \Z\setminus\{0\}\}.$$ 
If we study the attacking space of a single Queen placed at $Q\in \C$, then we write $A(Q)$ instead of $A(\{Q\})$. 

Convention: A Queen at square $Q$ does not attack $Q$, but it covers $Q$; i.e. if a Queen {\em attacks} $c-1$ squares, then it {\em covers} $c$ squares. For any given configuration $\C$, let $$\cover_n(\C)=\#\left([\C\cup A(\C)]\cap\B_n\right)$$ denote the total number of covered squares on $\B_n$, without counting any multiplicities of attacks. Hence, this is the usual cover that we wish to maximize, given a game board. 

A configuration $\C$ is $n$-feasible if $\C\subset \B_n$. We will only consider feasible Queen configurations.

A $q$-configuration $\C$ is $n$-{\em optimal}, if, for all $q$-configurations $\C'$, $\cover_n(\C')\le \cover_n(\C)$. Note that a $q$-configuration $\C$ does not depend on $n$ (provided feasible), but the question about its optimality does. If the board size is given we might drop the `$n$' and simply write ``optimal configuration''.  

A configuration $\C$ is {\em non-attacking} if it is pairwise non-attacking, that is, if, for all $\{Q, Q'\}\subset \C$, $Q'\not \in A(Q)$, which ny using our notation is equivalent to $\C \cap A(\C) = \varnothing$. 

The {\em border} of $\B_n$, $n\ge 2$, is the set $b_n=\B_n\setminus \B_{n-2}$. For example $b_2=\B_2$ and $$b_3=\{(-1,-1),(-1,0),(-1,1),(0,-1),(0,1),(1,-1),(1,0),(1,1)\}$$ In general, of course, $|\B_n|=n^2$ and $|b_n|= 4n-4$.

Non-trivial Queen configurations cannot be perfect in the sense that some squares will inevitably be attacked by more than one Queen. We begin by formalizing the number of times each square is attacked.

\begin{defi}[Attacking Number]
Consider a configuration $\C$. The attacking number of $s\in \Z\times \Z$ is $a_\C(s)=\#\{Q\in \C\mid s\in A(Q)\}$. 
\end{defi}




In order to quantify the inefficiency where there are some squares that are attacked by more than one Queen, we will define the concept of \emph{internal loss} of a configuration, normalized by the board size. We will use the convention that the attacking numbers smaller than two do not contribute to any {\em internal loss}. 

\begin{defi}[Internal Loss]\label{def:internalloss}
Consider a given configuration $\C\subset \B_n$. Its internal loss is $$\inl_n(\C)=\sum_{s\in A(\C)\cap\B_n} (a_\C(s)-1).$$
\end{defi}

For example, consider a $10\times 10$ board with Queen configuration, $\mathcal{C} = \{(-1,0), (0,2), (1,-1), (2,1)\}$.  The attacking numbers and internal loss are as in Figure~\ref{fig:loss}.

\begin{figure}[ht!]

\definecolor{attack0}{HTML}{FFFFFF}
\definecolor{attack1}{HTML}{DCE6F2}
\definecolor{attack2}{HTML}{B7C9E2}
\definecolor{attack3}{HTML}{7A9AC6}
\definecolor{attack4}{HTML}{3E5F99}

\centering

\begin{tikzpicture}[scale=0.55]

\foreach \x in {-5,...,4}{
\foreach \y in {-5,...,4}{

\pgfmathtruncatemacro{\attacks}{
((\x==0) || (\y==-2) || (abs(\x-0)==abs(\y+2))) +
((\x==-2) || (\y==-1) || (abs(\x+2)==abs(\y+1))) +
((\x==1) || (\y==0) || (abs(\x-1)==abs(\y-0))) +
((\x==-1) || (\y==1) || (abs(\x+1)==abs(\y-1)))
}

\ifcase\attacks
\fill[attack0] (\x,\y) rectangle ++(1,1);
\or
\fill[attack1] (\x,\y) rectangle ++(1,1);
\or
\fill[attack2] (\x,\y) rectangle ++(1,1);
\or
\fill[attack3] (\x,\y) rectangle ++(1,1);
\or
\fill[attack4] (\x,\y) rectangle ++(1,1);
\fi

}}

\draw[step=1cm,black,thin] (-5,-5) grid (5,5);

\foreach \x/\y in {0/-2,-2/-1,1/0,-1/1}{
\fill[attack0] (\x,\y) rectangle ++(1,1);
\node at (\x+0.5,\y+0.5) {\Large $\symqueen$};
}

\end{tikzpicture}
\hspace{1cm}
\begin{tikzpicture}

\draw (-0.2,-2.2) rectangle (3.2,1.8);

\node[align=center,font=\bfseries] at (1.5,1.3) {\textls[-60]Number of\\Queens attacking};

\node[left] at (0.3,0.4) {0};
\draw[fill=attack0,draw=black] (0.3,0.1) rectangle +(0.6,0.6);

\node[left] at (0.3,-0.4) {1};
\draw[fill=attack1,draw=black] (0.3,-0.7) rectangle +(0.6,0.6);

\node[left] at (0.3,-1.2) {2};
\draw[fill=attack2,draw=black] (0.3,-1.5) rectangle +(0.6,0.6);

\node[left] at (2,0.4) {3};
\draw[fill=attack3,draw=black] (2.1,0.1) rectangle +(0.6,0.6);

\node[left] at (2,-0.4) {4};
\draw[fill=attack4,draw=black] (2.1,-0.7) rectangle +(0.6,0.6);

\end{tikzpicture}

\caption{The attacking numbers that contribute to the internal loss $(2-1)\times 28 + (3-1)\times 4+(4-1)\times4 = 48$ of the Queen configuration in Figure~\ref{fig:4Qs}.}
\label{fig:loss}
\end{figure}


Consider a fixed number of Queens, $q$. A sequence of $q$-configurations $(\C_n)$ maps each $n$ to a $q$-configuration, where $\C_n\subset \B_n$. A sequence of $q$-configurations has {\em bounded internal loss} if there is a constant $\mu$ such that, for all $n$, $\inl(\C_n)\le \mu$. A fixed configuration $C$ has {\em bounded internal loss}, if the sequence $(\C_n)$ has bounded internal loss and where, for all $n$, $\C_n=\C$.

We will prove in Lemma~\ref{lem:nonattackboundedinloss} that a fixed Queen configuration has bounded internal loss if and only if it is a non-attacking configuration. 
 Hence, we will show that the internal loss for non-attacking configurations does not depend on $n$, for modestly large $n$.  In this study, we will usually assume that the board size is sufficiently large, so that the internal loss does not depend on $n$. If this is the case we may drop the ``$n$'' and write simply $\inl(\C)$.

In the following lemma and throughout the paper, we use the Chebyshev metric to measure distances on the board. The Chebyshev metric measures the maximum difference in any coordinate between two points; $d((x,y),(z,w))=\max\{|x-z|,|y-w|\}$. For example the Chebyshev distance of a point $(x,y)\in \mathbb{Z}^2$ to the central square $(0,0)$ is $d(x,y)=\max\{|x|,|y|\}$.

\begin{figure}[ht!]
\centering

\begin{tikzpicture}[scale = 0.55]
\fill[yellow!50] (3,4) rectangle (4,5);
\fill[red!50] (5, 10) rectangle (6, 11);

\foreach \x in {1,2,3,4, 5, 6}
    \draw[black!80] (\x,2) -- (\x,12);

\foreach \y in {3,4, 5, 6, 7, 8, 9, 10, 11}
    \draw[black!80] (0,\y) -- (7,\y);

\node at (3.5, 1.5) {0};
\node at (-0.5, 4.5) {0};
\node[circle, draw, fill = yellow!20, inner sep=2.5pt] at (3.5, 4.5) {};
\node[circle, draw, fill = red!20, inner sep=2.5pt] at (5.5, 10.5) {};
\draw[blue, line width=0.7pt] (3.5, 4.5) -- (3.5, 3.5) -- (5.5, 3.5);
\draw[blue, line width=0.7pt] (3.5, 4.5) -- (1.5, 4.5) -- (1.5, 6.5);
\draw[black, line width=0.7pt] (1.5, 6.5) -- (5.5, 3.5);
\node at (4.5, 3.2) {\textcolor{blue}{\small $\alpha$}};
\node at (1.2, 5.5) {\textcolor{blue}{\small $\alpha$}};
\node at (3.7, 5.35) {\small $2\alpha$};

\draw[red, line width=0.7pt] (1.5, 6.5) -- (5.4, 10.5) -- (5.4, 3.5);
\draw[green, line width=0.7pt] (3.5, 4.5) -- (5.6, 4.5) -- (5.6, 10.5);

\node at (1.5, 6.5) {\Large ${\symqueen}$};
\node at (5.5, 3.5) {\Large ${\symqueen}$};

\node at (5.6, 6.5) {-};
\node at (5.6, 4.5) {-};
\node at (5.6, 10.5) {-};

\node at (5.9, 5.5) {\textcolor{black}{\small $\alpha$}};
\node at (6, 8.5) {\textcolor{black}{\small $2\alpha$}};
\end{tikzpicture}
\caption{A non-attacking pair of Queens at distance $\alpha$ from the center and $2\alpha$ apart. Their attacking lines intersect at distance $3\alpha$ from the center.}
\label{fig:nonattackingboundedloss}
\end{figure}

\begin{lem}[Non-attacking Bounded Internal Loss]\label{lem:nonattackboundedinloss}
A fixed $q$-configuration $\C$ has bounded internal loss if and only if it is non-attacking.
\end{lem}
\begin{proof}
We will first consider a fixed configuration $\C$ that  contains a pair of attacking Queens.
As the board size increases, the common attacking line between this Queen pair extends across any newly added border. Consequently, the number of doubly covered squares grows without bound, and hence the internal loss of $\C$ is unbounded.

For the converse, suppose that $\C$ is non-attacking. Since $\C$ is fixed, all Queens lie within some finite distance $\alpha$ of the origin, where $\alpha = \max_{(x,y)\in \C}d(x,y)$. 
 Consider the boundary case where two extremal Queens lie at distance $\alpha$ from the origin, with internal distance $2\alpha$ between them. The diagonal attacking line from one Queen intersects an orthogonal  line from the other Queen at a distance at most $3\alpha$, as shown in Figure~\ref{fig:nonattackingboundedloss}. Therefore if we set the `radius' $\rho = 3\alpha+1$, then for every $(x,y)\in \B_n$ with $d(x,y)>\rho$, we have $a_\C(x,y)\le 1$. 
\end{proof}

Sometimes it is useful to distinguish between a Queen configuration, which assumes a given placement on some game board, and the observed (internal) pattern of the Queens. The term ``configuration'' refers to a given placement of the Queens on a fixed board, but the term ``pattern'' refers to the internal distances of the Queens, with the specific placement not fixed.

Next, let us study how the positioning of a single Queen on the board  influences its coverage. In particular, a Queen placed closer to the center of the board typically covers more squares. 
\begin{obs}\label{obs:central}
    A central single Queen on $\B_n$  covers $n+3(n-1)=4n-3$ squares, if $n$ is odd, and otherwise, if $n$ is even, she covers $n+2(n-1)+n-2=4n-4$ squares, since one diagonal is shorter.
\end{obs}


 Intuitively, as a Queen is placed further from the center, the number of covered squares decreases. To quantify this reduction, we define the notion of {\em centralized loss} by using the Chebyshev metric.

\begin{defi}[Centralized Loss]\label{def:centralloss}
Consider a single Queen $Q\in \B_n$. Its centralized loss, $\cl(Q)$, is zero (one) if $Q$ is centrally placed on a  board of odd (even) size. The centralized loss increases by two for each unit increase of Chebyshev-distance from a center square as seen in Figure~\ref{fig:centralizedLoss}. 
Let $\cl(\C)=\sum_{Q\in \C} \cl(Q)$ be the centralized loss of a configuration $\C$. 
\end{defi}
A fixed Queen configuration is {\em centralized} if it minimizes the centralized loss with respect to any  translation of the given Queen pattern. See Section~\ref{sec:app} for computational evidence and visualizations.

\begin{figure}[ht!]
\centering
\begin{tikzpicture}[scale=0.36]


\draw[gray,thin] (-2.5,-2.5) rectangle (2.5,2.5);

\foreach \x in {-1.5,-0.5,0.5,1.5}
  \draw[gray,thin] (\x,-2.5) -- (\x,2.5);
\foreach \y in {-1.5,-0.5,0.5,1.5}
  \draw[gray,thin] (-2.5,\y) -- (2.5,\y);

\foreach \x in {0}
  \foreach \y in {0}
   \node at (\x,\y) {\textcolor{gray}{\scriptsize 0}}; 

\foreach \x/\y in {
  -1/-1, -1/0, -1/1,
   0/-1,        0/1,
   1/-1,  1/0,  1/1}
  \node at (\x,\y) {\scriptsize 2};
  
\foreach \x/\y in {
-2/2,  -1/2,  0/2, 1/2,  2/2,
-2/1,                    2/1,
-2/0,                    2/0,
-2/-1,                   2/-1,
-2/-2, -1/-2, 0/-2, 1/-2, 2/-2}
  \node at (\x,\y) {\textbf{\scriptsize 4}}; 
\end{tikzpicture}
\hspace{5mm}
\begin{tikzpicture}[scale=0.36]


\draw[gray,thin] (-2.5,-2.5) rectangle (3.5,3.5);

\foreach \x in {-1.5, -0.5,0.5,1.5, 2.5}
  \draw[gray,thin] (\x,-2.5) -- (\x,3.5);
\foreach \y in {-1.5,-0.5,0.5,1.5, 2.5}
  \draw[gray,thin] (-2.5,\y) -- (3.5,\y);

\foreach \x/\y in {
0/1, 1/1,
0/0, 1/0}
  \node at (\x,\y) {\textcolor{gray}{\scriptsize 1}};

\foreach \x/\y in {
-1/2, 0/2, 1/2,    2/2,
-1/1,              2/1,
-1/0,              2/0,
-1/-1, 0/-1, 1/-1, 2/-1}
  \node at (\x,\y) {\scriptsize 3};
  
\foreach \x/\y in {
-2/3,  -1/3,  0/3, 1/3,  2/3, 3/3,
-2/2,                         3/2,
-2/1,                         3/1,
-2/0,                         3/0,
-2/-1,                        3/-1,
-2/-2, -1/-2, 0/-2,1/-2, 2/-2,3/-2}
  \node at (\x,\y) {\textbf{\scriptsize 5}};

\end{tikzpicture}






\caption{The centralized loss of a single Queen on the boards $\B_5$ and $\B_6$. }
\label{fig:centralizedLoss}
\end{figure}

\begin{lem}\label{lem:coverloss}
    A single Queen $Q\in\B_n$ satisfies $\cover_n(Q)=(4n-3) - \cl_n(Q)$.
\end{lem}

\begin{proof}
This is immediate by combining Observation~\ref{obs:central} with Definition~\ref{def:centralloss}. Namely, by inreasing a Queen's central distance with one unit, her cover decreases by one unit on each diagonal, but remains the same on both orthogonal attacking lines.
\end{proof}


Let us combine the two notions of ``loss'' introduced in Definitions~\ref{def:internalloss} and \ref{def:centralloss}. Intuitively, no other type of loss with respect to total cover is possible, as demonstrated in Theorem~\ref{thm:losscover}.  

\begin{defi}[Loss]\label{def:loss}
Let $\loss_n(\C) = \inl_n(\C) + \cl_n(\C)$ be the loss of a configuration $\C\subset\B_n$. 
\end{defi}
Sometimes we call ``loss'' by the more explicit {\em total loss}. Let us express the cover of a given configuration in terms of its loss. If the board is small relative to the configuration, then the loss may decrease, but the result implicitly covers those cases as well.\footnote{It is well-known that loss can decrease dramatically and allow for optimal configurations with many attacking Queen pairs, in terms of dominating the full board \cite{JW2004, W2002, TG2007, BM2002}.}
\begin{thm}\label{thm:losscover}
    For any fixed $q$-configuration $\C$, the coverage in $\B_n$ is $$\cover_n(\C)=(4n-3)q-\loss_n(\C).$$
\end{thm}

\begin{proof}
By Lemma~\ref{lem:coverloss}, a coverage measure that disregards any overlap reduction, is $(4n - 3)q - \cl(\C)$. But each square attacked by multiple Queens contributes only once to coverage, and the difference is given by $\inl(\C)$. Subtracting this overlap surplus yields the formula. 
\end{proof}
This in itself is a strong result, because of the exact correspondence between the cover of a configuration and its loss, which has generic complexity independent of $n$. But the theoretical value in our context will be more explicit in the following relaxation. 
\begin{cor}\label{cor:minlossmaxcov}
 Consider two $q$-configurations $\C$ and $\C'$ in $\B_n$. Then $\loss(\C)\ge \loss(\C')$ if and only if $\cover_n(\C)\le \cover_n(\C')$. In particular, any $q$-configuration that minimizes the loss on $\B_n$ maximizes the  coverage. 
\end{cor}
\begin{proof}
These are immediate consequences of Theorem~\ref{thm:losscover}.
\end{proof}

\section{The non-attacking threshold}\label{sec:nonattack}



In this section we study, for a fixed number of Queens $q$, properties of configurations that minimize the loss, which then, by Corollary~\ref{cor:minlossmaxcov},  maximize the coverage. We will consider sequences of Queen configurations, and unless otherwise stated, we regard a fixed configuration as a sequence of configurations on increasing board sizes. 
By using only the concept of centralized loss, we see that any single Queen of fixed location covers $4n-O(1)$ squares on the board $\B_n$. The contributions from all $q$ Queens, involve doubly and triply covered squares, and so on, and this induces further loss, relative to the amount $4nq$, in terms of the defined internal loss. 

A sequence of $q$-configurations has unbounded loss if there is no constant $k$ such that for all sufficiently large $n,  \loss(C_n)\le k$. 


\begin{lem}\label{lem:unbounded}
If a sequence of $q$-configurations, say $(\C_n)$, has unbounded internal loss on the corresponding sequence of boards $(\B_n)$, then, for all but finitely many board sizes, any fixed non-attacking $q$-configuration has greater cover.
\end{lem}
\begin{proof}
By Theorem~\ref{thm:losscover}, the cover of a $q$-configuration $\C$ on $\B_n$ is given by $$\cover_n(\C)=(4n-3)q-\loss(\C)=(4n-3)q-\inl(\C) - \cl(\C).$$ 

Let $\C$ be any fixed non-attacking $q$-configuration. By Lemma~\ref{lem:nonattackboundedinloss}, its internal loss is bounded, and since $\C$ is fixed, its centralized loss is constant. Hence $\loss(\C)$ is bounded, and so $\cover_n(\C) = 4qn - O(1)$.

Now let $(\C_n)$ be a sequence of $q$-configurations with unbounded internal loss. Then $\loss(\C_n)$ grows without bound and eventually exceeds any fixed constant. Hence, there must exist some  $N$ such that for all $n\ge N$,
$\cover_n(\C) > \cover_n(\C_n)$.
\end{proof}



Thus, a loss function that grows faster than a constant is a problem. It is possible to impose a bounded internal loss to sequences of configurations with infinitely many attacking configurations. To understand the case where internal loss is bounded, it is useful to define different types of positions within a sequence of configurations. A Queen in a sequence of configurations $\C_n$ is \emph{center-inclined} (\emph{corner-inclined}) if its Chebyshev distance to the central square (any corner of $\B_n)$ is bounded. Note that if the attacking Queens common attacking line is parallel to an edge of the board, regardless of placement, the internal loss will not be bounded as it guaranteed to increase as a function of $n$. 

\begin{lem}\label{lem:attackincorner}
Consider a sequence $(\C_n)$ that contains infinitely many $q$-configurations with attacking Queens. 
\begin{itemize}
    \item[(a)] If this sequence has a bounded internal loss, then for large $n$, the distance between a corner of $\B_n$ and any point on the line connecting any pair of attacking Queens must be bounded.
    \item[(b)]Any such sequence eventually has strictly smaller cover than any fixed non-attacking $q$-configuration.
\end{itemize}
\end{lem}

\begin{proof}
Let $(\C_n)$ be a sequence of $q$-configurations that have infinitely many configurations with attacking Queens.

(a) suppose that $(\C_n)$ has infinitely many configurations with a pair of attacking Queens such that the distance between any point on the line connecting them and a corner of $\B_n$ is not bounded. This implies that the length of the line is not bounded and the line size increases with $n$. In this case, such an attacking pair will double cover a number of squares that grows with $n$, and hence the internal loss is unbounded; see Figure~\ref{fig:attackQueencontdistance} (Left). 

(b) In the proof of Lemma~\ref{lem:unbounded} we saw that any fixed non-attacking $q$-configuration have $\cover_n(\C) = 4qn - O(1)$. On the other hand a sequence $(\C_n)$ with attacking pairs of Queens within a constant distance from a corner, loses almost one diagonal, and its centralized loss increases with $n$; see Figure~\ref{fig:attackQueencontdistance} (Right). Therefore, $\loss(\C_n)$ grows without bound and eventually exceeds any fixed constant. Hence, there must exist some  $N$ such that for all $n\ge N$,
$$\cover_n(\C) > \cover_n(\C_n).$$

\end{proof}

\begin{figure}[H]
\centering

\begin{minipage}{0.4\textwidth}
\centering
\begin{tikzpicture}[scale=1]
\fill[yellow!100] (1.85,1.85) rectangle (2.15, 2.15);
\node[circle, draw, fill = yellow!50, inner sep=1.5pt] at (2, 2) {};

\draw[thick] (0, 0) rectangle (4, 4);
\draw[dashed, thick] (-1, -1) rectangle (5, 5);
\draw[red, thick] (0, 0.5) -- (3.5,4);
\draw[dashed, thick, red] (0, 0.5) -- (-1, -0.5);
\draw[dashed, thick, red] (4.5, 5) -- (3.5, 4);

\node at (1,1.5) {\large $\symqueen$};
\node at (2.5,3) {\large $\symqueen$};

\node at (2, -0.5) {$\vdots$};
\node at (2, 4.5) {$\vdots$};
\node at (-0.5, 2) {$\cdots$};
\node at (4.5,2) {$\cdots$};

\end{tikzpicture}
\end{minipage}
\begin{minipage}{0.4\textwidth}
\centering
\begin{tikzpicture}

\fill[yellow!30] (1.85,2.85) rectangle (2.15,3.15);
\node[circle, draw, fill = yellow!10, inner sep=1.5pt] at (2, 3) {};
\fill[yellow!100] (2.35,2.35) rectangle (2.65,2.65);
\node[circle, draw, fill = yellow!50, inner sep=1.5pt] at (2.5, 2.5) {};
\draw[thick, dashed, yellow!70!orange] (2, 3) -- (2.5, 2.5);
\draw[thick] (0,0) rectangle (5,5);
\draw[thick, dashed] (0, 5) -- (-1, 6);
\draw[thick, dashed] (0, 0) -- (-1, 0) -- (-1, 6);
\draw[thick, dashed] (-1, 6) -- (5, 6) -- (5, 5);
\fill[red!30] (2.5,0) -- (5,2.5) -- (5, 0);
\draw[red, thick, <->] (3.1,0) -- (5,1.9);
\draw[thick, yellow!70!orange] (2.5, 2.5) -- (2.5, 0.8) -- (3.9,0.8);
\draw[thick, dashed, orange] (2, 3) -- (2, 0.8) -- (3.8,0.8);
\node at (3.9,0.8) {\large $\symqueen$};
\node at (4.7,1.6) {\large $\symqueen$};
\end{tikzpicture}
\end{minipage}
\caption{The pictures depict a pair of diagonally attacking Queens from a sequence of $q$-configurations. To the left: the board size increases, and the distances  of the two Queens with respect to the corners remains unbounded. Hence the internal loss grows linearly with $n$. To the right: the board size increases and the center is shifted in the top-left direction relative to the Queens. A pair of corner-inclined Queens stay within the bounding triangle (shaded red); their internal loss remains bounded, but the centralized loss increases.}
\label{fig:attackQueencontdistance}
\end{figure}

Let us state and prove our first main result. 
\begin{thm}[Non-attacking Threshold]\label{thm:main1}
For any given number of Queens $q$, there exists an $N_1(q)$ such that all optimal configurations on $\B_n$, for $n\ge N_1(q)$, have non-attacking Queens.
\end{thm}
\begin{proof}
For  a fixed $q$, as in the statement, consider a sequence of  optimal Queen configurations $(\C_{n})$. 
Suppose that infinitely many such configurations contain attacking pairs of Queens. By Lemma~\ref{lem:attackincorner}, any such sub-sequence eventually has strictly smaller cover than any fixed non-attacking $q$-configuration. This contradicts the optimality of all the configurations $\C_n$ for sufficiently large $n$. 
Therefore, there exists $N_1(q)$ such that for all $n \ge N_1(q)$, every optimal configuration on $\B_n$ is non-attacking.
\end{proof}
The {\em non-attacking threshold} is a number $N_1(q)$, such that every optimal $q$-configuration is non-attacking on $\B_{n}$ if $n\ge N_1(q)$. 

\section{The stabilizing threshold}\label{sec:stablize}





Experimental results yield explicit non-attacking thresholds for $q=2,3,4$, with $N_1(q)=10,9,10$ respectively. However, in general, a non-attacking threshold does not imply stabilization of optimal configurations. Here, stabilization means that the set of optimal configurations remains the same as the board size increases.  Experimentally for $q=2,3,4$, and all $n\ge 10, 12, 15$ respectively, sets of fixed configurations have been established via exhaustive search. Motivated by methods developed in this section, by placing Queens in a central area while computing coverages on sufficiently large boards, we have been able to verify further fixed sets of configurations, for $5\le q\le 9$. See Section~\ref{sec:app} for computational results and visualizations. 
   
We show that optimal configurations stabilize and become constant, but clearly this could not happen before the non-attacking threshold from Theorem~\ref{thm:main1} has settled. Hence, in this section, we exclusively study non-attacking configurations, and we prove the existence of a stabilizing threshold, $N_2(q)$, for any $q\ge 2$. Recall Figure~\ref{fig:4Qs}, where $q=4$. This is an example for which $N_1(q)< N_2(q)$; the two upper left pictures show the only optimal configurations for the $10\times10$ board (modulo symmetry), but these configurations are not optimal for the $12\times12$ board (see two lower left pictures). Conversely, the configuration shown in the lower right picture is optimal for the $12\times12$ board but is not optimal for the $10\times10$ board.


\begin{thm}[Stabilizing Threshold]\label{thm:main2}
Consider any $q \ge 2$. Then there exists a threshold $N_2(q)$ such that,  for all even $n \ge N_2(q)$ (and analogously for odd $n$), the set $\mathcal O_n$ of optimal $q$-configurations on $\B_n$ is constant.
\end{thm}

\begin{proof}
For a $q$-configuration $\C$, let $\rho(\C)$ denote the Chebyshev distance from the center of $\B_n$ to the farthest Queen in $\C$. By Corollary~\ref{cor:minlossmaxcov}, it suffices to study the minimum loss; every configuration that minimizes the loss maximizes the coverage. We must  establish that the set of configurations with minimum loss remains constant, for sufficiently large board sizes, depending only on the number of Queens, $q$, and the parity of the board size.\\

\noindent {\em Claim}: There is a finite maximum radius, $\rho_{\max,q}$, such that 
\[
\rho_{\max,q} := \max \{ \rho(\C) \mid \C \text{ is a $q$-configuration with minimum loss} \}.
\]

Before we prove this claim, let us explain how it implies the result. 
Given this maximum radius, we derive a threshold $N_2(q)$. In analogy with the proof of Lemma~\ref{lem:nonattackboundedinloss}, a pair of Queens placed at distances no more than $\rho$ from the center can have common attacking squares up to an additional distance of $2\rho$. Thus, all crossings between attacking lines are confined to a central square of side length at most $6\rho$. To ensure that these interactions remain within the board $\B_n$, it suffices to require that
\[
n \ge 6\rho_{max, q} + 1.
\]
Define $N_2(q) = 6\rho_{max, q} + 1$. Then, for all even $n \ge N_2(q)$ (and analogously for odd $n$), the loss function is fully stabilized and unaffected by the boundary of $\B_n$. Since the set of $q$-configurations of minimal loss is finite and $n$-independent for $n \ge N_2(q)$, the set of optimal configurations $\mathcal O_n$ remains constant for all such $n$.\\

\noindent {\em Proof of Claim}: 
For each fixed $q$, there exists a configuration with finite loss such as the stairs configuration $\S_q$ (for example see 
Remark ~\ref{rem:Nq} just below). Thus any loss minimizing configuration $\C_{\min}$ must satisfy
\begin{equation}\label{eq:minloss}
    \loss(\C_{\min}) \le \loss(\S_q).
\end{equation}
Since the centralized loss of any $\C_{\min}$ contributes to $\loss(\C_{\min})$, inequality~\eqref{eq:minloss} bounds $\rho(\C_{\min})\le \cl(\C_{\min})\le  \loss(\S_q)$. Thus, all loss-minimizing $q$-configurations are inside a fixed finite region of the board. Hence $\rho_{\mathrm{maxmin}}(q)$ is finite and well defined. 
Further, since the centralized loss differs between odd and even boards, the set $\mathcal O_n$ of optimal configurations, as well as the value of $\rho_{q,\max} $, can depend on the parity of the board size.
\end{proof}


Here, we abuse the notation $N_2(q)$ to mean any number that suffices as a threshold (whenever we require uniqueness we would obviously take the smallest such number).

We believe that the stairs construction from the proof is near optimal for large boards, so let us look into this construction a bit and study its loss components empirically. 

\begin{figure}[H]
    \centering
    \includegraphics[width=1\linewidth]{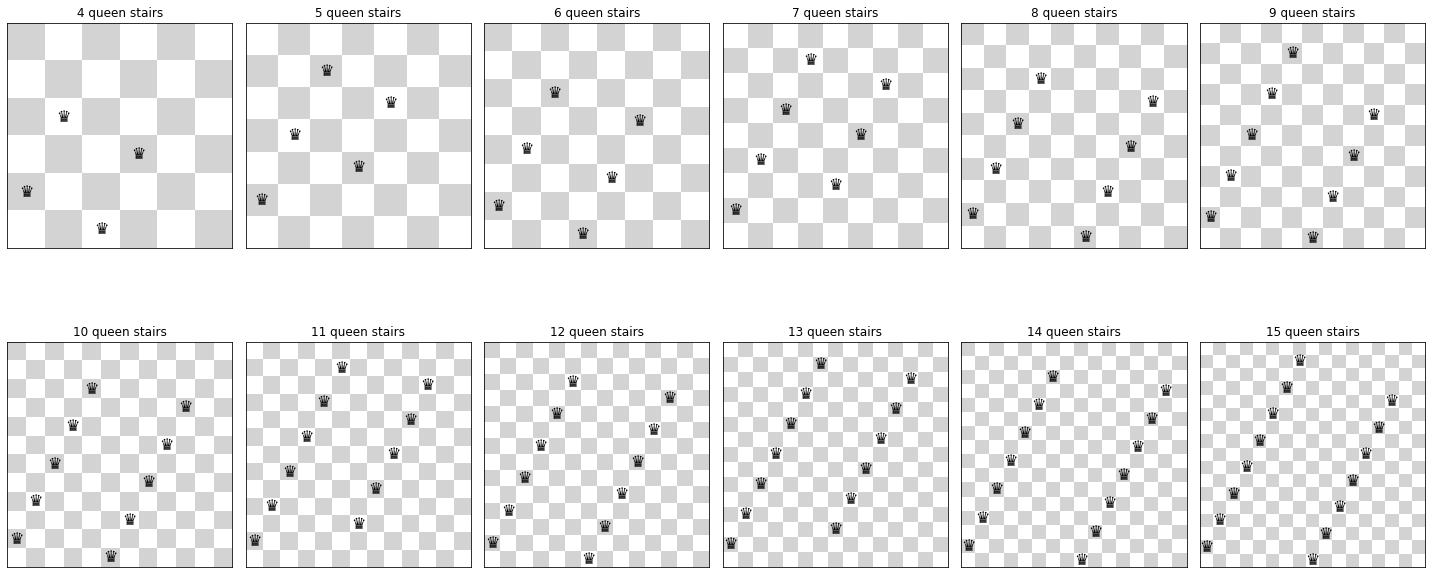}
    \caption{Some $q$-stairs configurations.}
    \label{fig:qstairs}
\end{figure}

\begin{figure}[H]
\includegraphics[width=.5\linewidth]{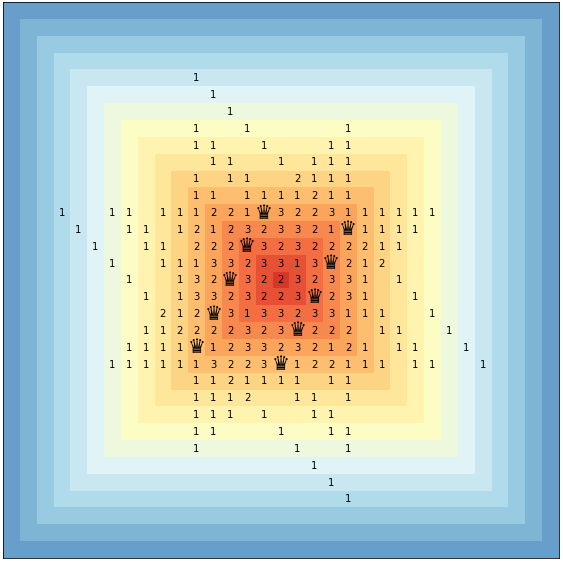}
\caption{A stairs construction with ten Queens. Here centrality loss (in colors) sums up to $70$, while the internal loss contributes almost five times as much, $340$. See also Table~\ref{tab:stair}.}
\label{fig:inlosscenloss}
\end{figure}

\begin{remark}\label{rem:Nq}
Empirically, many constructions satisfy $2\rho_{q,\max} \le  q+a $, for small $a\ge 1$. A classical configuration is a ``$q$-stairs-construction'', $\S_q$. It consists of two sequences of Queens that evolve with Knight's jumps. There are a few possible definitions, with small variations in cover efficiency depending on $n$ and $q$, but a couple of things are in common. The distance between the sequences varies depending on $q$ modulo~$6$ and some of these Queen-stair pairs fit inside $\B_q$ (the $q$ by $q$ board), while other ones require the somewhat larger $\R_q$ (the $q$ by $q+1$ board; see Section~\ref{sec:app}). See Figure~\ref{fig:qstairs} for an example, with corresponding losses presented in Table~\ref{tab:stair}; see also Figure~\ref{fig:inlosscenloss}, which illustrates centrality loss together with the internal loss in a stairs configuration. 
This suggests that $N_2(q)$ may be chosen close to $3q + 4$ in many cases. 
\end{remark}

\begin{table}[H]\caption{Loss table for a $q$-stairs construction, for even and odd large board sizes, respectively. The parity of the board size matters for the centrality loss, because there are four central $1$-loss squares for even $n$, but a single $0$-loss square, for odd $n$.
}\label{tab:stair}
\centering
\begin{tabular}{c|c||c|c|}
$q$ & Internal loss & Centrality odd & Total loss odd\\
\hline
2 & 10 & 4 & 14 \\
3 & 27 & 8 & 35 \\
4 & 48 & 12 & 60 \\
5 & 76 & 16 & 92 \\
6 & 116 & 26 & 142 \\
7 & 158 & 32 & 190 \\
8 & 222 & 50 & 272 \\
9 & 277 & 60 & 337 \\
10 & 340 & 70 & 410 \\
11 & 410 & 80 & 490 \\
12 & 496 & 100 & 596 \\
13 & 580 & 112 & 692 \\
14 & 698 & 144 & 842 \\
15 & 791 & 160 & 951 \\
16 & 896 & 176 & 1072 \\
\hline
\end{tabular}
\begin{tabular}{|c|c|}
 Centrality even & Total loss even\\
\hline
  4 & 14 \\
  7 & 34 \\
 12 & 60 \\
 17 & 93 \\
 26 & 142 \\
 33 & 191 \\
 50 & 272 \\
 59 & 336 \\
 70 & 410 \\
 81 & 491 \\
 100 & 596 \\
 113 & 693 \\
 144 & 842 \\
 159 & 950 \\
 176 & 1072 \\
 \hline
\end{tabular}
\end{table}

Let us return to think about the proof of Theorem ~\ref{thm:main2}; it does not rely on the staircase configuration $\S_q$ being optimal or even close to optimal. The only requirement is that, for each fixed $q$, it has finite loss that becomes independent of $n$ once the board is sufficiently large. This gives an upper bound on the loss of any optimal configuration: no configuration with loss greater than $\loss(\S_q)$ can be optimal. 
By this section, all optimal configurations lie inside a ``fixed box around the center of the board'', depending only on the number of Queens, $q$; the set of optimal configurations stabilizes for all $n \ge N_2(q)$. In the coming sections, we shed some more light on this fenomenon. In particular, we begin by disection the notion of internal loss.

\section{Balance versus concentration}\label{sec:balovercon}
We will continue to discuss how both types of loss tend to attract Queens into a small rectangle $\R_q$ of size $q\times(q+1)$, supporting Remark~\ref{rem:Nq}. A defined balanced non‑attacking configuration can always be placed within this region. While such attraction is natural for centralized loss, it is perhaps less obvious that internal conflicts lead to the same effect.

Recall that there is an internal loss due to crossing attack lines (Definition~\ref{def:internalloss}). Such internal loss emanating from a single pair of Queens is upper bounded by $12$ squares. 
An interesting observation is as follows. Let us denote by $\Delta(x_1,x_2)=x_1-x_2\pmod 2$. We call a Queen $Q$ {\em even} if $\Delta(Q)=0$ and otherwise {\em odd}. A pair of Queens is {\em congruent} if the Queens have the same parity. 

\begin{obs}
For all sufficiently large boards, a single pair of non-attacking Queens has internal loss either $10$, if they are non-congruent (i.e. if $\Delta(Q_1)\neq\Delta(Q_2)\pmod 2$), or $12$ (otherwise).
\end{obs}
Due to this observation, we are interested in maximizing the number of non-congruent pairs of Queens, as to minimize the induced internal loss with respect to crossings of attack lines. The recurrence defined by $a_q+a_{q-1}=\binom{q}{2}$, with $a(1)=0$, is the same as the sequence of ``quarter squares'', for $q>0$, $a_q=\lfloor\frac{q^2}{4}\rfloor$:
$$0,1,2,4,6,9,12,16,20,25,30,36,42,\ldots$$

We call a configuration of $q$ Queens {\em balanced} if the difference between the number of odd and even Queens is at most one. If $q=3$ we can have at most two pairs of non-congruent Queens. If $q=4$, we get at most four pairs, and if $q=5$, we get at most six such pairs. Indeed, the multiplication principle applies.

\begin{prop}\label{lem:parity}
The sequence $(a_q)$, $q=1,2,\ldots$, represents the maximum number of non-congruent pairs of Queens in a $q$-configuration, that is the number of  $(Q_1,Q_2)\in\C_q\times \C_q$ such that $\Delta(Q_1)\ne \Delta(Q_2)$. 
\end{prop}
\begin{proof} 
Any balanced configuration maximizes the number of non-congruent Queen pairs. Suppose there are $e$ even Queens and $o$ odd ones. The number of congruent pairs is $\binom{e}{2}+\binom{o}{2}$ and the number of non-congruent pairs is $e\cdot o$. The latter expression is maximized when $e=o$ if $q$ is even, and $e=o-1$ or $e=o+1$ if $q$ is odd. 

For the second part, we use the identity $\binom{q+1}{2}-\binom{q-1}{2}=q-1$. The number of congruent pairs increase by $q-1$ because the odd pairs increase by one and the even pairs increase by one, i.e. if $o+e=q-1$ then $o+1+e+1=q+1$, and   $\binom{e+1}{2}+\binom{o+1}{2} - \binom{e}{2}-\binom{o}{2}=q-1$. Similarly, we get $(o+1)(e+1)-oe=q$, so the number of non-congruent pairs increases by $q=\binom{q+2}{2}-\binom{q}{2}$. Thus, the sequences coincide.
\end{proof}
However there is a trade-off to the degree of `balance' that we will next discuss. In particular, there are optimal configurations in stabilized sets that are unbalanced. 
Internal loss at a square can reach as high as three, as demonstrated already in Figure~\ref{fig:loss} for four Queens. This connects to a somewhat counterintuitive feature of Queen placements. 

For any fixed numbers of even and odd parity Queens, \(e\) and \(o\) respectively, the constant  
\begin{align}\label{eq:gamma}
    \gamma = 12\binom{e}{2} + 12\binom{o}{2} + 10eo
\end{align}  
represents the total number of crossings of Queen pairs, counting multiplicities. 
In a given configuration, multiple pairs may overlap on the same square. We can express $\gamma$ as 
\(
\sum_{s \in A} \binom{a(s)}{2},
\) 
where \(a(s)\) is the number of Queens attacking square \(s\). For instance, if four Queens attack a single square, that square accounts for `six overlapping pairs'. 

This motivates the notion of a discrepancy in terms of \emph{overlap concentration}, defined as, for a given configuration, 
\[
\eta := \sum_{s \in A} \left[\binom{a(s)}{2}-(a(s) - 1)\right]. 
\]
 
Note that when a square is attacked by \(a(s)\) Queens, it contributes \(a(s) - 1\) to internal loss while \(\binom{a(s)}{2}\) represents the number of pairs that meet at square $s$.  As \(\binom{a(s)}{2}\) grows quadratically while \(a(s) - 1\) grows linearly, larger attacking numbers make the overlap more efficient in this sense, by increasing \(\eta\) as it `uses up pairs'. In contrast, if every overlapping pair would correspond to a distinct square, then \(\eta\) approaches its hypothetical minimum value, which is zero. 

Although \(\binom{a}{2}\) grows quadratically with \(a\), the attacking number is bounded due to the non-attacking constraint: no square can be attacked by more than four Queens, for otherwise there is an attacking pair of Queens. The difference of \(\binom{a}{2}\) relative to $a-1$ still favors higher attacking numbers, as $a=2,3,4$ correspond to the discrepancy term $0, 1$ and $3$, respectively. 

Thus, concentrating overlaps at fewer high-attacking-number squares tends to reduce internal loss relative to the total number of overlapping interactions. Naively this might appear as a structural advantage, but there is an apparent trade-off in that producing many high attack crossing may reduce balance, which then would increase the internal loss. A well centralized balanced configuration still seems like a good condidate for optimality.


To summarize, we have the following proposition. 

\begin{prop}\label{prop:inlgammaeta}
For any configuration $\C$, we have $\inl(\C) = \gamma(\C) - \eta(\C)$.  
\end{prop}
\begin{proof}
This is immediate by the definitions of these quantities and by noting that $\gamma = \sum_{s \in A} \binom{a(s)}{2}$.
\end{proof}

Observe that an immediate consequence of this result is that, if we have two configurations $\C$ and $\C'$ with the same inloss, then $\gamma(\C') - \gamma(\C) = \eta(\C') - \eta(\C)$. Thus, to preserve the internal loss (or equivalently the loss if centralization is the same), any deviation in balance must be compensated precisely in terms of overlap concentration. 

\begin{cor} 
Consider $q$-configurations $\C$ and $\C'$ such that $\inl(\C) = \inl(\C')$. Then any difference in balance is compensated precisely by the configurations difference in overlap concentration, that is $\gamma(\C') - \gamma(\C) = \eta(\C') - \eta(\C)$.
\end{cor} 
\begin{proof}
    This is immediate by Proposition~\ref{prop:inlgammaeta}.
\end{proof}

To illustrate the effects of parity balance versus overlap concentration, Figure~\ref{fig:5Queenloss} displays the attacking numbers for the squares that contribute to internal loss in two optimal five-Queen configurations. While the configuration on the right is unbalanced in parity, it compensates through a larger overlap concentration.


\begin{figure} [ht!]
\begin{center}
\begin{tikzpicture} [scale=0.5]
\draw (-1.5,0.5) node {$\symqueen$};
\draw (1.5,-0.5) node {$ \symqueen$};
\draw (-0.5,-1.5) node {$ \symqueen$};
\draw (0.5,1.5) node {$ \symqueen$};
\draw (2.5,2.5) node {$ \symqueen$};
\draw[step=1cm,gray, thin] (-5,-5) grid (7,7);
\draw[color=blue] (-1.5,-3.5) node {$ 2 $};
\draw[color=blue] (2.5,-3.5) node {$ 2 $};
\draw[color=blue] (-1.5,-2.5) node {$ 2 $};
\draw[color=blue] (1.5,-2.5) node {$ 2 $};
\draw[color=red] (-1.5,-1.5) node {$ 3 $};
\draw[color=blue] (-3.5,-1.5) node {$ 2 $};
\draw[color=ForestGreen] (0.5,-1.5) node {$ 4 $};
\draw[color=blue] (-2.5,0.5) node {$ 2 $};
\draw[color=ForestGreen] (-1.5,-0.5) node {$ 4 $};
\draw[color=blue] (-2.5,-0.5) node {$ 2 $};
\draw[color=ForestGreen] (-0.5,-0.5) node {$ 4 $};
\draw[color=red] (-0.5,0.5) node {$ 3 $};
\draw[color=blue] (-1.5,1.5) node {$ 2 $};
\draw[color=ForestGreen] (-0.5,1.5) node {$ 4 $};
\draw[color=blue] (-2.5,1.5) node {$ 2 $};
\draw[color=ForestGreen] (0.5,0.5) node {$ 4 $};
\draw[color=red] (-1.5,2.5) node {$ 3 $};
\draw[color=red] (0.5,2.5) node {$ 3 $};
\draw[color=blue] (-3.5,2.5) node {$ 2 $};
\draw[color=ForestGreen] (1.5,0.5) node {$ 4 $};
\draw[color=blue] (-1.5,3.5) node {$ 2 $};
\draw[color=red] (1.5,3.5) node {$ 3 $};
\draw[color=red] (2.5,0.5) node {$ 3 $};
\draw[color=blue] (2.5,4.5) node {$ 2 $};
\draw[color=blue] (4.5,0.5) node {$ 2 $};
\draw[color=blue] (-1.5,6.5) node {$ 2 $};
\draw[color=blue] (1.5,-3.5) node {$ 2 $};
\draw[color=blue] (-0.5,-2.5) node {$ 2 $};
\draw[color=red] (0.5,-0.5) node {$ 3 $};
\draw[color=blue] (1.5,-1.5) node {$ 2 $};
\draw[color=red] (2.5,-1.5) node {$ 3 $};
\draw[color=red] (2.5,-0.5) node {$ 3 $};
\draw[color=red] (1.5,1.5) node {$ 3 $};
\draw[color=red] (3.5,1.5) node {$ 3 $};
\draw[color=red] (1.5,2.5) node {$ 3 $};
\draw[color=blue] (4.5,2.5) node {$ 2 $};
\draw[color=blue] (5.5,-0.5) node {$ 2 $};
\draw[color=blue] (2.5,-4.5) node {$ 2 $};
\draw[color=blue] (-2.5,-1.5) node {$ 2 $};
\draw[color=blue] (0.5,-2.5) node {$ 2 $};
\draw[color=red] (2.5,1.5) node {$ 3 $};
\draw[color=blue] (-3.5,1.5) node {$ 2 $};
\draw[color=blue] (3.5,-1.5) node {$ 2 $};
\draw[color=red] (-0.5,2.5) node {$ 3 $};
\draw[color=blue] (3.5,2.5) node {$ 2 $};
\draw[color=blue] (-4.5,2.5) node {$ 2 $};
\draw[color=blue] (6.5,-1.5) node {$ 2 $};
\draw[color=blue] (-0.5,5.5) node {$ 2 $};
\draw[color=blue] (2.5,3.5) node {$ 2 $};
\draw[color=blue] (0.5,4.5) node {$ 2 $};
\end{tikzpicture}\hspace{.5cm}
\begin{tikzpicture} [scale=0.5]
\draw (-1.5,1.5) node {$ \symqueen$};
\draw (2.5,-0.5) node {$ \symqueen$};
\draw (-0.5,-1.5) node {$ \symqueen$};
\draw (1.5,2.5) node {$ \symqueen$};
\draw (0.5,0.5) node {$ \symqueen$};
\draw[step=1cm,gray, thin] (-5,-5) grid (6,6);
\draw[color=blue] (-1.5,-4.5) node {$ 2 $};
\draw[color=blue] (-1.5,-2.5) node {$ 2 $};
\draw[color=blue] (2.5,-2.5) node {$ 2 $};
\draw[color=red] (-1.5,-1.5) node {$ 3 $};
\draw[color=blue] (-4.5,-1.5) node {$ 2 $};
\draw[color=ForestGreen] (1.5,-1.5) node {$ 4 $};
\draw[color=blue] (-3.5,1.5) node {$ 2 $};
\draw[color=ForestGreen] (-1.5,-0.5) node {$ 4 $};
\draw[color=blue] (-3.5,-0.5) node {$ 2 $};
\draw[color=ForestGreen] (0.5,-0.5) node {$ 4 $};
\draw[color=blue] (-1.5,0.5) node {$ 2 $};
\draw[color=red] (-2.5,0.5) node {$ 3 $};
\draw[color=ForestGreen] (-0.5,0.5) node {$ 4 $};
\draw[color=red] (-0.5,1.5) node {$ 3 $};
\draw[color=red] (-1.5,2.5) node {$ 3 $};
\draw[color=ForestGreen] (-0.5,2.5) node {$ 4 $};
\draw[color=blue] (-2.5,2.5) node {$ 2 $};
\draw[color=ForestGreen] (0.5,1.5) node {$ 4 $};
\draw[color=blue] (-1.5,3.5) node {$ 2 $};
\draw[color=red] (0.5,3.5) node {$ 3 $};
\draw[color=red] (1.5,1.5) node {$ 3 $};
\draw[color=blue] (1.5,4.5) node {$ 2 $};
\draw[color=ForestGreen] (2.5,1.5) node {$ 4 $};
\draw[color=blue] (-1.5,5.5) node {$ 2 $};
\draw[color=blue] (2.5,5.5) node {$ 2 $};
\draw[color=blue] (4.5,1.5) node {$ 2 $};
\draw[color=blue] (2.5,-4.5) node {$ 2 $};
\draw[color=red] (-0.5,-0.5) node {$ 3 $};
\draw[color=blue] (-0.5,-3.5) node {$ 2 $};
\draw[color=red] (0.5,-2.5) node {$ 3 $};
\draw[color=red] (1.5,-0.5) node {$ 3 $};
\draw[color=red] (2.5,-1.5) node {$ 3 $};
\draw[color=blue] (3.5,-1.5) node {$ 2 $};
\draw[color=blue] (2.5,0.5) node {$ 2 $};
\draw[color=red] (3.5,0.5) node {$ 3 $};
\draw[color=ForestGreen] (1.5,0.5) node {$ 4 $};
\draw[color=blue] (4.5,-0.5) node {$ 2 $};
\draw[color=red] (2.5,2.5) node {$ 3 $};
\draw[color=blue] (5.5,2.5) node {$ 2 $};
\draw[color=blue] (2.5,3.5) node {$ 2 $};
\draw[color=blue] (-2.5,-1.5) node {$ 2 $};
\draw[color=blue] (1.5,-3.5) node {$ 2 $};
\draw[color=blue] (0.5,-1.5) node {$ 2 $};
\draw[color=blue] (3.5,2.5) node {$ 2 $};
\draw[color=blue] (-4.5,2.5) node {$ 2 $};
\draw[color=blue] (5.5,-1.5) node {$ 2 $};
\draw[color=blue] (-0.5,4.5) node {$ 2 $};
\draw[color=blue] (0.5,2.5) node {$ 2 $};
\end{tikzpicture}
\end{center}
\caption{Whenever centralized, both these Queen-patterns are optimal on large boards. To the right we see a tradeoff while being un-balanced, this is compensated by a larger proportion of high attacking numbers. Inloss: to the left, $30+2\times 14+3\times 6$; to the right: $28+2\times 12+3\times 8$.}
\label{fig:5Queenloss}
\end{figure}

\begin{remark}
    Since the balanced condition can be attained inside $\R_q$ (a $q$ by $q+1$ rectangle), we guess that centralized loss can be trivially bounded as $2q(q+1)$. Thus heuristically, the total cover is no less than $\eta = 4nq-\gamma-2q(q+1)$.  Suppose $n=3q+3$, and $e=o= q/2$. Then $$\eta = 4(3q+3)q-2q(q+1)- 6q(q/2-1)-5q/2= 7q^2+13.5q$$ If $n=3q+4$, then we get instead the cover $\eta=7q^2+17.5q$. Consider 8 Queens on a 27 by 27 board. We get $\eta = 7\times 64+13.5\times 8= 556$. Compare with stairs construction, $(4n-3)\times 8-264=(4\times 27-3)\times 8-264=576$.
\end{remark}

As mentioned, there is a theoretical worst-case scenario for internal loss, analogous to the ideal best-case scenario for centralized loss. Since high attacking numbers are advantageous, the worst-case internal loss arises when all overlapping attacks involve exactly two Queens. In this setting, each such pair contributes one unit of internal loss, and the total internal loss reaches the bound $\gamma$ (for a fixed number of Queens of respective parity) as defined in~\eqref{eq:gamma}. 

\begin{figure}[ht!]
\includegraphics[width=.45\linewidth]{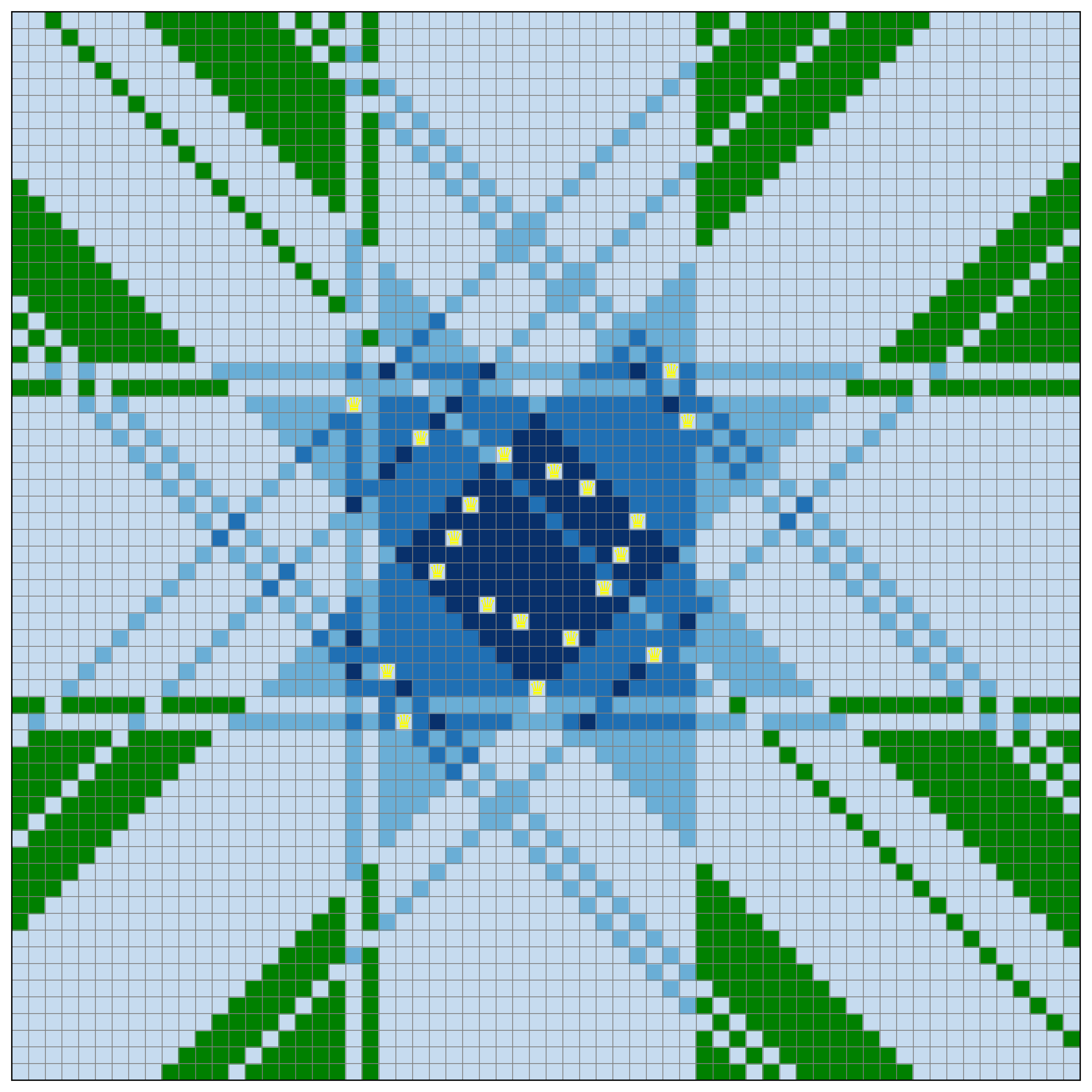}
\caption{An illustration of attacking numbers for a configuration of twenty Queens on $\B_{65}$. Here green color indicates non-covered areas, and various shades of blue indicate positive attacking numbers, with dark blue representing $a(
s)=4$.}\label{fig:20Q}
\end{figure}

Let us give some further intuition along the lines in this section: attacking concentration cannot be increased by moving Queens outwards on a large board. Conversly one could compress any large pattern with fixed internal loss towards the centre, modulo discretization, to minimize centralized loss. But, because of the discreteness of the problem, it is not obvious how far you might be able to push a certain pattern simulataneously towards the centre. A better understanding of this interdependency might lead to reasonable explicit bounds of the thresholds. Does this phenomenon help to maximize the number of Queen pairs with multiple (near four) interactions? In Figure~\ref{fig:4Qs} we see how the number of squares with four interactions remein constant, while the four Queens have been moved towards the center. For large number of Queens, we do not yet have any explicit understanding of the various optimizing Queen patterns on large boards. In Figure~\ref{fig:20Q}, we see a reasonably good configuration, produced by one of our many semi-optimizing codes, and we note that the most central squares have been ignored by the Queens. Is this a common feature for a many-Queen-configuration to be optimal on large board sizes? High attacking numbers are concentrated in the central part, which points at an empirical tendency; inloss and centralized loss both benefit by concentrating Queens near the middle of the board, perhaps skipping a relatively small number of central squares. 

\section{Computations}\label{sec:app}

Solutions of the classical non-attacking $q$-Queens problem, are all subsets of $\B_q$. However, already for $q=2,3$, there is no non-attacking $q$-configuration in $\B_q$. Moreover, we have several examples where a classical configuration may be slightly shifted while remaining cover-optimal. Most interestingly, we find that the fundamental classical solution for $q=6$ is not optimal for large board sizes, but the only optimal solutions require a slightly larger area. Thus, we search for non-attacking solutions on a somewhat larger board, namely the $q\times(q+1)$-rectangle $\R_q=\{-q/2+1,\ldots , q/2\}\times \{-q/2+1,\ldots , q/2+1\}$, if $q$ is even, and otherwise $\R_q=\{-(q-1)/2,\ldots , (q-1)/2-1\}\times \{-(q-1)/2,\ldots , (q-1)/2\}$. Sometimes we cannot find balanced (classical) $q$-configurations inside $\B_q$. However, by extending to $\R_q$, this is always possible, and by our experiments, for $q\le 9$ there is always a balanced optimal solution. Such optimal solutions may sometimes be shifted, while still remaining centralized (see in particular case $q=6$).


Let us illustrate the idea behind Theorem~\ref{thm:main2}, by finding the set of {\em $q$-cover-fundamental solutions} for $q=2,3,\ldots, 9$, respectively. Those are the representatives modulo symmetry of the corresponding equivalence classes; those appear as the Dihedral group, or subgroups thereof, since, if we find one solution, we may always reflect or rotate it to obtain another one with the same cover (some of the configurations have inherent symmetries, which then reduces to some subgroup by $2$ or $4$). 

We will see that some of those configurations include  classical ``fundamental solutions'' on $\B_q$ (the $q\times q$ board). But there are interesting exceptions, for example for $q=6$. For visibility, we draw a red square around the relevant part of the picture, respectively. We draw a full game board, of size the experimental lower bound of the stabilizing threshold. For $q=2,3,4$, we have run exhaustive search (in C code) to obtain a repetition of the configurations for at least two even and odd board sizes beyond the threshold, respectively. 
For larger $q$, exhaustive search over all $n^2 \choose q$ configurations becomes intractable. We therefore reformulate the problem as an integer linear program. For each square $v \in B_n$, we introduce a binary decision variable $x_v$ indicating whether a Queen is placed at $v$, together with auxiliary variables capturing the attacking number at each square and its contribution to the internal loss. The objective is to minimize the loss function from Definition~\ref{def:loss}, subject to the cardinality constraint $\sum_v x_v = q$. We solve this model using Gurobi's solution pool feature, which enumerates all configurations attaining the optimal loss value. This allows us to rigorously certify the sets of cover-fundamental solutions reported in this section for $2 \leq q \leq 9$. Empirically, we observe that the stabilizing threshold $N_2(q)$ is about $3q+3$, or a bit larger in some cases.

For $q=2$, $n\ge 10=N_2(2)$, there are two cover-fundamental solutions for even $n$, each containing 8 symmetric configurations, and (just below) odd $n$ is similar.\\

\begin{center}
\chessboardtwoQueen{}{\node at (-0.5,-0.5) {{\tiny \symqueen}}; \node at (.5,1.5) {{\tiny \symqueen}};}
\chessboardtwoQueen{}{\node at (0.5,-0.5) {{\tiny \symqueen}}; \node at (1.5,1.5) {{\tiny \symqueen}};}\\
\end{center}

\begin{center}
\chessboardthreeQueenodd{}{\node at (-0.5,-0.5) {{\tiny \symqueen}}; \node at (.5,1.5) {{\tiny \symqueen}};}
\chessboardthreeQueenodd{}{\node at (-0.5,-1.5) {{\tiny \symqueen}}; \node at (.5,.5) {{\tiny \symqueen}};}\\
\end{center}

\vspace{5mm}
For $q=3$, and even $n\ge 12$ there is only one representative, with 8 symmetric configurations:\\

\begin{center}
\chessboardthreeQueen{}{\node at (0.5,-0.5) {{\tiny \symqueen}}; \node at (1.5,1.5) {{\tiny \symqueen}};\node at (-1.5,.5) {{\tiny \symqueen}};}\\
\end{center}

\vspace{5mm}
For $q=3$, and odd $n\ge 13$ there are four representatives, each one with 8 symmetric configurations:\\

\begin{center}
\chessboardthreeQueenodd{}{\node at (0.5,0.5) {{\tiny \symqueen}}; \node at (-1.5,-.5) {{\tiny \symqueen}};\node at (1.5,-1.5) {{\tiny \symqueen}};}
\chessboardthreeQueenodd{}{\node at (0.5,-.5) {{\tiny \symqueen}}; \node at (1.5,1.5) {{\tiny \symqueen}};\node at (-1.5,.5) {{\tiny \symqueen}};}
\chessboardthreeQueenodd{}{\node at (-0.5,1.5) {{\tiny \symqueen}}; \node at (1.5,.5) {{\tiny \symqueen}};\node at (-1.5,-.5) {{\tiny \symqueen}};}
\chessboardthreeQueenodd{}{\node at (0.5,-0.5) {{\tiny \symqueen}}; \node at (-.5,1.5) {{\tiny \symqueen}};\node at (-1.5,-1.5) {{\tiny \symqueen}};}\\
\end{center}
\vspace{5mm}
For $q=4$, with $n\ge 15$ even,  there are two representatives, one central with two (mirror symmetric) configurations, and another one with 8 configurations, and if $n$ odd there is only one class, with $8$ symmetric configurations. Let us illustrate the more interesting case for even $n$, with a representative for each class:\\

\begin{center}
\chessboardfourQueeneven{}{\node at (0.5,1.5) {{\tiny \symqueen}}; \node at (-1.5,.5) {{\tiny \symqueen}};\node at (1.5,-.5) {{\tiny \symqueen}}; \node at (-.5,-1.5) {{\tiny \symqueen}};}
\chessboardfourQueeneven{}{\node at (0.5,.5) {{\tiny \symqueen}}; \node at (-1.5,-.5) {{\tiny \symqueen}};\node at (1.5,-1.5) {{\tiny \symqueen}}; \node at (-.5,-2.5) {{\tiny \symqueen}};}\\
\end{center}

\vspace{15mm}

For $q=5$, and even $n\ge  18$ there are 5  representatives, with 8 in each class (total 40):\\
\begin{center}
\chessboardfiveQueen{}{\node at (0.5,1.5) {{\tiny \symqueen}}; \node at (-1.5,.5) {{\tiny \symqueen}};\node at (1.5,-.5) {{\tiny \symqueen}}; \node at (-.5,-1.5) {{\tiny \symqueen}}; \node at (-2.5,-2.5) {{\tiny \symqueen}};}
\chessboardfiveQueen{}{\node at (1.5,1.5) {{\tiny \symqueen}}; \node at (-.5,.5) {{\tiny \symqueen}};\node at (2.5,-.5) {{\tiny \symqueen}}; \node at (.5,-1.5) {{\tiny \symqueen}}; \node at (-1.5,-2.5) {{\tiny \symqueen}};}
\chessboardfiveQueen{}{\node at (0.5,2.5) {{\tiny \symqueen}}; \node at (-1.5,1.5) {{\tiny \symqueen}};\node at (1.5,.5) {{\tiny \symqueen}}; \node at (-.5,-.5) {{\tiny \symqueen}}; \node at (-2.5,-1.5) {{\tiny \symqueen}};}
\chessboardfiveQueen{}{\node at (1.5,2.5) {{\tiny \symqueen}}; \node at (-.5,1.5) {{\tiny \symqueen}};\node at (2.5,.5) {{\tiny \symqueen}}; \node at (.5,-.5) {{\tiny \symqueen}}; \node at (-1.5,-1.5) {{\tiny \symqueen}};}\vspace{1.5mm}
\chessboardfiveQueen{}{\node at (.5,-2.5) {{\tiny \symqueen}}; \node at (-1.5,1.5) {{\tiny \symqueen}};\node at (1.5,.5) {{\tiny \symqueen}}; \node at (-.5,-.5) {{\tiny \symqueen}}; \node at (-2.5,-1.5) {{\tiny \symqueen}};}
\\
\end{center}
\vspace{10 mm}
For $q=5$, and odd $n\ge  17$ there are 2 representatives, with 8 and 2 elements respectively (total 10):\\
\begin{center}
\chessboardfiveQueenodd{}{\node at (1.5,2.5) {{\tiny \symqueen}}; \node at (-.5,1.5) {{\tiny \symqueen}};\node at (2.5,.5) {{\tiny \symqueen}}; \node at (.5,-.5) {{\tiny \symqueen}}; \node at (-1.5,-1.5) {{\tiny \symqueen}};}
\chessboardfiveQueenodd{}{\node at (1.5,-1.5) {{\tiny \symqueen}}; \node at (-.5,2.5) {{\tiny \symqueen}};\node at (2.5,1.5) {{\tiny \symqueen}}; \node at (.5,.5) {{\tiny \symqueen}}; \node at (-1.5,-.5) {{\tiny \symqueen}};}
\\
\end{center}
For $q=6$, and odd $n\ge 21$ there are four representatives (covering 346 squares for $\B_{21}$), each with 8 elements:\\
\begin{center}
\chessboardsixQueen{}{ \node at (-2.5,-2.5) {{\tiny \symqueen}};\node at (-.5,.5) {{\tiny \symqueen}};\node at (-1.5,3.5) {{\tiny \symqueen}};\node at (.5,2.5) {{\tiny \symqueen}};\node at (2.5,1.5) {{\tiny \symqueen}};\node at (1.5,-.5) {{\tiny \symqueen}};}
\chessboardsixQueen{}{ \node at (-2.5,-2.5) {{\tiny \symqueen}};\node at (-.5,1.5) {{\tiny \symqueen}};\node at (-1.5,3.5) {{\tiny \symqueen}};\node at (1.5,2.5) {{\tiny \symqueen}};\node at (2.5,.5) {{\tiny \symqueen}};\node at (.5,-.5) {{\tiny \symqueen}};}\\
\chessboardsixQueen{}{ \node at (-2.5,-3.5) {{\tiny \symqueen}};\node at (-.5,-.5) {{\tiny \symqueen}};\node at (-1.5,2.5) {{\tiny \symqueen}};\node at (.5,1.5) {{\tiny \symqueen}};\node at (2.5,.5) {{\tiny \symqueen}};\node at (1.5,-1.5) {{\tiny \symqueen}};}
\chessboardsixQueen{}{ \node at (-2.5,-3.5) {{\tiny \symqueen}};\node at (-.5,.5) {{\tiny \symqueen}};\node at (-1.5,2.5) {{\tiny \symqueen}};\node at (1.5,1.5) {{\tiny \symqueen}};\node at (2.5,-.5) {{\tiny \symqueen}};\node at (.5,-1.5) {{\tiny \symqueen}};}
\\
\end{center}

Similarly, for $q=6$, and even $n\ge  22$ there are four representatives (covering 370 squares for $\B_{22}$), each with 8 elements:\\
\begin{center}
\chessboardsixQueeneven{}{ \node at (-2.5,-3.5) {{\tiny \symqueen}};\node at (-.5,-.5) {{\tiny \symqueen}};\node at (-1.5,2.5) {{\tiny \symqueen}};\node at (.5,1.5) {{\tiny \symqueen}};\node at (2.5,.5) {{\tiny \symqueen}};\node at (1.5,-1.5) {{\tiny \symqueen}};}
\chessboardsixQueeneven{}{ \node at (-2.5,-3.5) {{\tiny \symqueen}};\node at (-.5,.5) {{\tiny \symqueen}};\node at (-1.5,2.5) {{\tiny \symqueen}};\node at (1.5,1.5) {{\tiny \symqueen}};\node at (2.5,-.5) {{\tiny \symqueen}};\node at (.5,-1.5) {{\tiny \symqueen}};}\\
\chessboardsixQueeneven{}{ \node at (-3.5,-3.5) {{\tiny \symqueen}};\node at (-1.5,-.5) {{\tiny \symqueen}};\node at (-2.5,2.5) {{\tiny \symqueen}};\node at (-.5,1.5) {{\tiny \symqueen}};\node at (1.5,.5) {{\tiny \symqueen}};\node at (.5,-1.5) {{\tiny \symqueen}};}
\chessboardsixQueeneven{}{ \node at (-3.5,-3.5) {{\tiny \symqueen}};\node at (-1.5,.5) {{\tiny \symqueen}};\node at (-2.5,2.5) {{\tiny \symqueen}};\node at (.5,1.5) {{\tiny \symqueen}};\node at (1.5,-.5) {{\tiny \symqueen}};\node at (-.5,-1.5) {{\tiny \symqueen}};}
\\
\end{center}
Observe that all cover-fundamental solutions for $q=6$ fits in a central rectangle of size $q\times (q+1)$, but cannot be fitted inside the classical board of size $q\times q$, with its different single fundamental solution.
\clearpage

For seven Queens, with odd board sizes $n\ge 25$, the red frame gets size 7 by 7, with $32=8+8+4+4+8$ elements. We get: 

\begin{center}
\chessboardsevenQueen{.17}{%
  \node at (-2.5,-2.5) {{\tiny \symqueen}};
  \node at (-.5,-1.5) {{\tiny \symqueen}};
  \node at (1.5, -.5) {{\tiny \symqueen}};
  \node at (3.5, .5) {{\tiny \symqueen}};
  \node at ( -1.5,1.5) {{\tiny \symqueen}};
  \node at ( .5,2.5) {{\tiny \symqueen}};
  \node at ( 2.5, 3.5) {{\tiny \symqueen}};}
\chessboardsevenQueen{.17}{%
  \node at (-2.5,-2.5) {{\tiny \symqueen}};
  \node at (2.5,-1.5) {{\tiny \symqueen}};
  \node at (0.5,-0.5) {{\tiny \symqueen}};
  \node at (-1.5,0.5) {{\tiny \symqueen}};
  \node at (3.5,1.5) {{\tiny \symqueen}};
  \node at (1.5,2.5) {{\tiny \symqueen}};
  \node at (-0.5,3.5) {{\tiny \symqueen}};
}
\chessboardsevenQueen{.17}{%
  \node at (-2.5,-1.5) {{\tiny \symqueen}};
  \node at (-1.5,1.5) {{\tiny \symqueen}};
  \node at (-0.5,-2.5) {{\tiny \symqueen}};
  \node at (0.5,0.5) {{\tiny \symqueen}};
  \node at (1.5,3.5) {{\tiny \symqueen}};
  \node at (2.5,-0.5) {{\tiny \symqueen}};
  \node at (3.5,2.5) {{\tiny \symqueen}};}
\end{center}
\begin{center}
\chessboardsevenQueen{.17}{%
  \node at (-2.5, -1.5) {{\tiny \symqueen}};
  \node at (-1.5,  1.5) {{\tiny \symqueen}};
  \node at (-0.5,  3.5) {{\tiny \symqueen}};
  \node at ( 0.5,  0.5) {{\tiny \symqueen}};
  \node at ( 1.5, -2.5) {{\tiny \symqueen}};
  \node at ( 2.5, -0.5) {{\tiny \symqueen}};
  \node at ( 3.5,  2.5) {{\tiny \symqueen}};
}
\chessboardsevenQueen{.17}{%
  \node at (-2.5, -1.5) {{\tiny \symqueen}};
  \node at (-1.5,  2.5) {{\tiny \symqueen}};
  \node at (-0.5, -0.5) {{\tiny \symqueen}};
  \node at ( 0.5,  3.5) {{\tiny \symqueen}};
  \node at ( 1.5,  0.5) {{\tiny \symqueen}};
  \node at ( 2.5, -2.5) {{\tiny \symqueen}};
  \node at ( 3.5,  1.5) {{\tiny \symqueen}};
}
\end{center}

Note, there is one fundamental configuration among the classical ones, that is non-optimal in our sense, that is, it is not cover-fundamental: 
$[(-3, -2), (-2, 0), (-1, -3), (0, 3), (1, 1), (2, -1), (3, 2)]$.

For seven Queens, with even board sizes $n\ge 24$, the red frame gets size $8$ by $8$, with $128=32+32+16+16+32$ elements. We get the same patterns as for odd board sizes, but with $4$ respective $2$ shifted placements of the respective patterns.\\

\begin{center}
\chessboardsevenQueeneven{.17}{%
  \node at (-2.5,-2.5) {{\tiny \symqueen}};
  \node at (-.5,-1.5) {{\tiny \symqueen}};
  \node at (1.5, -.5) {{\tiny \symqueen}};
  \node at (3.5, .5) {{\tiny \symqueen}};
  \node at ( -1.5,1.5) {{\tiny \symqueen}};
  \node at ( .5,2.5) {{\tiny \symqueen}};
  \node at ( 2.5, 3.5) {{\tiny \symqueen}};}
\chessboardsevenQueeneven{.17}{%
  \node at (-2.5,-2.5) {{\tiny \symqueen}};
  \node at (2.5,-1.5) {{\tiny \symqueen}};
  \node at (0.5,-0.5) {{\tiny \symqueen}};
  \node at (-1.5,0.5) {{\tiny \symqueen}};
  \node at (3.5,1.5) {{\tiny \symqueen}};
  \node at (1.5,2.5) {{\tiny \symqueen}};
  \node at (-0.5,3.5) {{\tiny \symqueen}};
}
\chessboardsevenQueeneven{.17}{%
  \node at (-2.5,-1.5) {{\tiny \symqueen}};
  \node at (-1.5,1.5) {{\tiny \symqueen}};
  \node at (-0.5,-2.5) {{\tiny \symqueen}};
  \node at (0.5,0.5) {{\tiny \symqueen}};
  \node at (1.5,3.5) {{\tiny \symqueen}};
  \node at (2.5,-0.5) {{\tiny \symqueen}};
  \node at (3.5,2.5) {{\tiny \symqueen}};}\\
\end{center}
\begin{center}
\chessboardsevenQueeneven{.17}{%
  \node at (-2.5, -1.5) {{\tiny \symqueen}};
  \node at (-1.5,  1.5) {{\tiny \symqueen}};
  \node at (-0.5,  3.5) {{\tiny \symqueen}};
  \node at ( 0.5,  0.5) {{\tiny \symqueen}};
  \node at ( 1.5, -2.5) {{\tiny \symqueen}};
  \node at ( 2.5, -0.5) {{\tiny \symqueen}};
  \node at ( 3.5,  2.5) {{\tiny \symqueen}};
}
\chessboardsevenQueeneven{.17}{%
  \node at (-2.5, -1.5) {{\tiny \symqueen}};
  \node at (-1.5,  2.5) {{\tiny \symqueen}};
  \node at (-0.5, -0.5) {{\tiny \symqueen}};
  \node at ( 0.5,  3.5) {{\tiny \symqueen}};
  \node at ( 1.5,  0.5) {{\tiny \symqueen}};
  \node at ( 2.5, -2.5) {{\tiny \symqueen}};
  \node at ( 3.5,  1.5) {{\tiny \symqueen}};
}
\end{center}

\newpage For eight Queens, with even board sizes, $n\ge 28$, there are six equivalence classes, each of size $8$, so in total $48$ optimal configurations. Odd board sizes have a smaller central bounding box, of size nine by nine, and otherwise the same patterns and sizes of equivalence classes. There are three cover-fundamental Queen patterns, and they all fit inside $\B_8$, while there are twelve classical fundamental non-attacking solutions. 

\begin{center}
\begin{tabular}{cc}
\drawchessboard{}{%
  \node at (-3.5,-4.5) {{\tiny \symqueen}}; \node at (-2.5,-0.5) {{\tiny \symqueen}}; \node at (-1.5,2.5) {{\tiny \symqueen}};
  \node at (-0.5,.5) {{\tiny \symqueen}}; \node at (0.5,-2.5) {{\tiny \symqueen}}; \node at (1.5,1.5) {{\tiny \symqueen}};
  \node at (2.5,-3.5) {{\tiny \symqueen}}; \node at (3.5,-1.5) {{\tiny \symqueen}};
}
&
\drawchessboard{}{%
  \node at (-3.5,-3.5) {{\tiny \symqueen}}; \node at (-2.5,.5) {{\tiny \symqueen}}; \node at (-1.5,3.5) {{\tiny \symqueen}};
  \node at (-0.5,1.5) {{\tiny \symqueen}}; \node at (0.5,-1.5) {{\tiny \symqueen}}; \node at (1.5,2.5) {{\tiny \symqueen}};
  \node at (2.5,-2.5) {{\tiny \symqueen}}; \node at (3.5,-.5) {{\tiny \symqueen}};
}\\
 
\drawchessboard{}{%
  \node at (-3.5,-4.5) {{\tiny \symqueen}}; \node at (-2.5,1.5) {{\tiny \symqueen}}; \node at (-1.5,-1.5) {{\tiny \symqueen}};
  \node at (-0.5,.5) {{\tiny \symqueen}}; \node at (0.5,2.5) {{\tiny \symqueen}}; \node at (1.5,-3.5) {{\tiny \symqueen}};
  \node at (2.5,-.5) {{\tiny \symqueen}}; \node at (3.5,-2.5) {{\tiny \symqueen}};
}
&
\drawchessboard{}{%
  \node at (-3.5,-3.5) {{\tiny \symqueen}}; \node at (-2.5,1.5) {{\tiny \symqueen}}; \node at (-1.5,3.5) {{\tiny \symqueen}};
  \node at (-0.5,-1.5) {{\tiny \symqueen}}; \node at (0.5,2.5) {{\tiny \symqueen}}; \node at (1.5,-.5) {{\tiny \symqueen}};
  \node at (2.5,-2.5) {{\tiny \symqueen}}; \node at (3.5,.5) {{\tiny \symqueen}};
}\\

\drawchessboard{}{%
  \node at (-3.5,-3.5) {{\tiny \symqueen}}; \node at (-2.5,1.5) {{\tiny \symqueen}}; \node at (-1.5,-2.5) {{\tiny \symqueen}};
  \node at (-0.5,0.5) {{\tiny \symqueen}}; \node at (0.5,2.5) {{\tiny \symqueen}}; \node at (1.5,-.5) {{\tiny \symqueen}};
  \node at (2.5,-4.5) {{\tiny \symqueen}}; \node at (3.5,-1.5) {{\tiny \symqueen}};
}
&
\drawchessboard{}{%
  \node at (-3.5,-2.5) {{\tiny \symqueen}}; \node at (-2.5,2.5) {{\tiny \symqueen}}; \node at (-1.5,-1.5) {{\tiny \symqueen}};
  \node at (-0.5,1.5) {{\tiny \symqueen}}; \node at (0.5,3.5) {{\tiny \symqueen}}; \node at (1.5,.5) {{\tiny \symqueen}};
  \node at (2.5,-3.5) {{\tiny \symqueen}}; \node at (3.5,-.5) {{\tiny \symqueen}};
}\\
\end{tabular}
\end{center}
\newpage
There are nine cover-fundamentals for the case $q=9$ and $n\ge 31$ odd. In total there are $64$ configurations, and exactly two classes have four elements each, while the remaining seven classes all have eight elements. All these solutions fit inside $\B_9$, and exactly one Queen pattern out of the nine is un-balanced.\\
\begin{center}
\chessboardnineQueenodd{}{%
  \node at (-3.5, -3.5) {\tiny \symqueen};
  \node at (-2.5, -0.5) {\tiny \symqueen};
  \node at (-1.5, 1.5) {\tiny \symqueen};
  \node at (-0.5, 3.5) {\tiny \symqueen};
  \node at (0.5, -2.5) {\tiny \symqueen};
  \node at (1.5, 0.5) {\tiny \symqueen};
  \node at (2.5, -1.5) {\tiny \symqueen};
  \node at (3.5, 4.5) {\tiny \symqueen};
  \node at (4.5, 2.5) {\tiny \symqueen};
} 
\chessboardnineQueenodd{}{%
  \node at (-3.5, -3.5) {\tiny \symqueen};
  \node at (-2.5, 0.5) {\tiny \symqueen};
  \node at (-1.5, 4.5) {\tiny \symqueen};
  \node at (-0.5, 1.5) {\tiny \symqueen};
  \node at (0.5, -0.5) {\tiny \symqueen};
  \node at (1.5, -2.5) {\tiny \symqueen};
  \node at (2.5, 3.5) {\tiny \symqueen};
  \node at (3.5, -1.5) {\tiny \symqueen};
  \node at (4.5, 2.5) {\tiny \symqueen};
} 
\chessboardnineQueenodd{}{%
  \node at (-3.5, -3.5) {\tiny \symqueen};
  \node at (-2.5, 1.5) {\tiny \symqueen};
  \node at (-1.5, 3.5) {\tiny \symqueen};
  \node at (-0.5, -1.5) {\tiny \symqueen};
  \node at (0.5, 2.5) {\tiny \symqueen};
  \node at (1.5, -0.5) {\tiny \symqueen};
  \node at (2.5, -2.5) {\tiny \symqueen};
  \node at (3.5, 4.5) {\tiny \symqueen};
  \node at (4.5, 0.5) {\tiny \symqueen};
}\\ 
\chessboardnineQueenodd{}{%
  \node at (-3.5, -3.5) {\tiny \symqueen};
  \node at (-2.5, 2.5) {\tiny \symqueen};
  \node at (-1.5, -0.5) {\tiny \symqueen};
  \node at (-0.5, 3.5) {\tiny \symqueen};
  \node at (0.5, -1.5) {\tiny \symqueen};
  \node at (1.5, 0.5) {\tiny \symqueen};
  \node at (2.5, 4.5) {\tiny \symqueen};
  \node at (3.5, -2.5) {\tiny \symqueen};
  \node at (4.5, 1.5) {\tiny \symqueen};
} 
\chessboardnineQueenodd{}{%
  \node at (-3.5, -2.5) {\tiny \symqueen};
  \node at (-2.5, -0.5) {\tiny \symqueen};
  \node at (-1.5, 4.5) {\tiny \symqueen};
  \node at (-0.5, 2.5) {\tiny \symqueen};
  \node at (0.5, 0.5) {\tiny \symqueen};
  \node at (1.5, -1.5) {\tiny \symqueen};
  \node at (2.5, -3.5) {\tiny \symqueen};
  \node at (3.5, 1.5) {\tiny \symqueen};
  \node at (4.5, 3.5) {\tiny \symqueen};
} 
\chessboardnineQueenodd{}{%
  \node at (-3.5, -2.5) {\tiny \symqueen};
  \node at (-2.5, 0.5) {\tiny \symqueen};
  \node at (-1.5, 2.5) {\tiny \symqueen};
  \node at (-0.5, 4.5) {\tiny \symqueen};
  \node at (0.5, -1.5) {\tiny \symqueen};
  \node at (1.5, 1.5) {\tiny \symqueen};
  \node at (2.5, -0.5) {\tiny \symqueen};
  \node at (3.5, -3.5) {\tiny \symqueen};
  \node at (4.5, 3.5) {\tiny \symqueen};
}\\ 
\chessboardnineQueenodd{}{%
  \node at (-3.5, -2.5) {\tiny \symqueen};
  \node at (-2.5, 1.5) {\tiny \symqueen};
  \node at (-1.5, -3.5) {\tiny \symqueen};
  \node at (-0.5, 2.5) {\tiny \symqueen};
  \node at (0.5, 0.5) {\tiny \symqueen};
  \node at (1.5, -1.5) {\tiny \symqueen};
  \node at (2.5, 4.5) {\tiny \symqueen};
  \node at (3.5, -0.5) {\tiny \symqueen};
  \node at (4.5, 3.5) {\tiny \symqueen};
} 
\chessboardnineQueenodd{}{%
  \node at (-3.5, -2.5) {\tiny \symqueen};
  \node at (-2.5, 1.5) {\tiny \symqueen};
  \node at (-1.5, 4.5) {\tiny \symqueen};
  \node at (-0.5, -1.5) {\tiny \symqueen};
  \node at (0.5, 0.5) {\tiny \symqueen};
  \node at (1.5, 3.5) {\tiny \symqueen};
  \node at (2.5, -0.5) {\tiny \symqueen};
  \node at (3.5, -3.5) {\tiny \symqueen};
  \node at (4.5, 2.5) {\tiny \symqueen};
} 
\chessboardnineQueenodd{}{%
  \node at (-3.5, -2.5) {\tiny \symqueen};
  \node at (-2.5, 3.5) {\tiny \symqueen};
  \node at (-1.5, 0.5) {\tiny \symqueen};
  \node at (-0.5, -1.5) {\tiny \symqueen};
  \node at (0.5, 4.5) {\tiny \symqueen};
  \node at (1.5, 1.5) {\tiny \symqueen};
  \node at (2.5, -0.5) {\tiny \symqueen};
  \node at (3.5, -3.5) {\tiny \symqueen};
  \node at (4.5, 2.5) {\tiny \symqueen};
} 
\end{center}
In the case of $n\ge 28$ even, there are instead $256$ configurations induced by the same nine patterns, while the size of the central bounding box increases to $\B_{10}$. The $256$ elements consist of eight goups of size four and $28$ of size eight. We omit the pictures, since they are immediate variations of the above nine ones.

\section{A historical Queen problem perspective}\label{sec:history}

The study of Queens on chessboards has followed two classical traditions that,
for much of their history, developed separately. One tradition concerns the
placement of mutually non-attacking Queens. Although the problem was already
studied privately by Gauss in the 1830s, it entered the public mathematical
literature in 1848 when Max Bezzel posed the famous problem of arranging eight
non-attacking Queens on the standard chessboard. Franz Nauck soon generalized
Bezzel’s problem to arbitrary $n\times n$ boards, giving rise to the
$n$-Queens problem. By the late nineteenth century, explicit pattern-based
constructions were known that produce at least one placement of $n$ mutually
non-attacking Queens on an $n\times n$ board for every $n \ge 4$
\cite{Nauck1850,MB1850,EL1883}. These classical results addressed feasibility
alone. With one Queen in each row and column, such configurations trivially
dominate the board, so questions of coverage efficiency, overlap of attack
lines, or optimality under limited resources do not arise. In this respect
they are analogous to Euler’s celebrated knight’s tour problem of 1759
\cite{Euler1759}, which likewise concerns the existence of a single global
structure rather than optimization or coverage.

A different tradition, emerging in the twentieth century and developed
systematically by Cockayne, Hedetniemi, and others, studied the Queen
domination problem: what is the minimum number of Queens required to dominate
an $n\times n$ board so that every square is occupied or attacked
\cite{CH1977,C1990}? In contrast to the classical ``$n$-Queens problem'', mutual
non-attack is not required, and allowing attacks among Queens can reduce the
number needed for domination. Subsequent constructive work, notably by Burger
and Mynhardt \cite{BM2002}, established explicit asymptotic upper bounds,
showing that at most $8n/15+O(1)$ Queens suffice for domination, while Weakley
and others computed exact domination numbers for small boards \cite{W2002}.
Related variants include domination by non-attacking Queens on infinite
boards \cite{DS2019}, as well as computational approaches based on
metaheuristics such as genetic algorithms \cite{AV2017}.

In contrast, our work studies small fixed sets of Queens tasked with covering
as much of an expanding board as possible. Surprisingly, this approach
naturally combines ideas from both earlier traditions. On small boards,
Queen configurations that include attacking pairs can be optimal. As the
board grows, such attacking placements become inefficient because fully
overlapping attack lines result in an internal loss that increases linearly
with the side length of the board. Thus, to keep the internal loss bounded,
Queens become non-attacking, giving rise to a first threshold. Once
non-attacking placements are achieved, certain Queen patterns are singled
out, and a stabilizing threshold appears. These stable Queen patterns remain
optimal, provided they are suitably centralized, as the board expands, and
interestingly they are often subsets of classical non-attacking solutions,
although there are notable exceptions. Our study shows how, when domination
of the entire board is relaxed to the task of maximizing coverage, the
classical focus on non-attacking placements and the later emphasis on
domination efficiency naturally join forces. For a fixed number of Queens
$q$, we study sets of densely packed, coverage-maximizing
non-attacking Queen configurations.

\section{Discussion and future work}\label{sec:fut}
It remains an open problem to establish reasonable bounds for the thresholds $N_1(q)$ and $N_2(q)$: can we determine exact values or tight asymptotics? For instance, does $N_1(q)$ grow linearly with $q$, or follow a sublinear or superlinear pattern? Even a rough estimate like $N_1(q) = \Theta(q)$ would be informative. 
In the case of the stabilizing threshold, the inherent tension between balance and overlap concentration seems especially relevant.

Problem: establish rigorously our stable sets for $2\le q\le 9$, and find more stabilizing sets, for $q>9$. 

There is a third threshold, $N_3$, such that, given a fixed number of Queens, on smaller boards, every optimal Queen configuration must have each Queen attacking at least one other. For instance, on a \( 3 \times 3 \) board, two non-attacking Queens can cover the entire board, so the threshold for two Queens must be the trivial \( 2 \times 2 \) case. 
Now consider \( 50 \) Queens on a \( 9 \times 9 \) board. Placing one Queen in the upper-left corner leaves a non-attacking region of size 
\[
2(7 + 6 + \cdots + 1) = 56 > 49,
\]
sufficient to place the remaining Queens disjointly while still covering the board. This trick fails for an \( 8 \times 8 \) board, suggesting that the threshold might be 8 in this case. 

This is another threshold, $N_4$, between the ones discussed in this paper. Given a number of Queens $q$, for what board size is there an optimal configuration that will remain optimal for all larger board sizes? Not that, for all $q$, $N_3(q)< N_1(q) \le N_4(q)\le N_2(q)$. We conjecture that $N_4(4)=11$.

We leave the general question of such thresholds open for further exploration.

In our problems, we fix the number of Queens; hence the search space remains polynomial in \( n \), following the formula
\[
\binom{n^2}{q} \sim n^{2q}.
\]
One could also study configurations in which the number of Queens grows proportionally with \( n \), say \( q = \alpha n \) for some \( \alpha > 0 \). In this regime, the configuration space grows super-polynomially, and new geometric or asymptotic phenomena may arise. It would be interesting to construct such sequences of Queen configurations that achieve optimal or near-optimal coverage as \( n \to \infty \). Observe that, if $\alpha <1/3$, then the heuristics from this paper imply that, for large $n$, every set of optimal configurations is the same as the `constant set' described here, although the problem is differently phrased.  

We can also extend this framework to other pieces like Rooks, Bishops, Knights, or bounded-range Queens, with vision limited by a radius or decreasing proportionally with distance. This would model realistic scenarios where influence range is physically constrained, such as in communication systems or security monitoring.

From an algorithmic angle, we can ask how many distinct optimal configurations exist for fixed $q$ and intermediate sized $n$ (say $q<n<3q$), and whether they can be enumerated efficiently. Experiments have shown that deterministic greedy strategies such as placing queens sequentially to cover the most new cells are often sub optimal, but nevertheless they could guide in proving approximation guarantees. Various probabilistic greedy algorithms can complement such efforts.

Finally, one can generalize from counting covered cells to maximizing a weighted coverage, where each cell has an importance weight $w(i,j)$. What configurations maximize total weighted coverage? This could model applications in surveillance or rescue systems where some areas are more critical than others.

Each of these directions extends the theory of directional influence over grids and raises rich questions at the interface of combinatorics, optimization, and spatial systems.

\section{Further applications}\label{sec:thoughts}

Our twist of the Queen cover problem has potential applications in several domains where directional influence over a grid is essential. In sensor networks, for example, Queens can model directional sensors whose coverage extends along rows, columns, and diagonals. Classical work on sensor placement and coverage optimization has emphasized the importance of maximizing area coverage with limited resources (\cite{Tossa2020, Njoya2017}). Studies of connectivity and barrier coverage in grid-based deployments further highlight parallels with our stabilization results \cite{LYC2015}. In structured environments like warehouses or data centers, deploying a fixed number of such sensors to maximize monitoring coverage is a key challenge. Our results help in estimating the maximum possible area that can be covered for a given number of sensors, and show that beyond a certain grid size, the optimal configurations naturally become non-overlapping, interference-free, and eventually stable, meaning the same placement remains optimal for all larger grids. This can simplify large-scale deployment planning.

A further application arises in game AI, where our model parallels the concept of \emph{influence maps} used in grid-based strategy games \cite{Tozour2002}. In such settings, each agent projects influence over the cells it can attack or observe, and optimal deployment aims to maximize territorial control or visibility. A Queen’s coverage pattern is a direct analogue: its row, column, and diagonal control correspond to an agent’s influence zone. Our stabilization results suggest that beyond a certain map size, optimal influence-maximizing placements converge to fixed patterns, implying that precomputed strategies may remain effective across larger arenas without recomputation.

Another application is in urban emergency response planning, such as placing directional fire alarms or beacons in a grid of rooms or intersections. Beyond a certain grid size, optimal placement patterns stabilize, allowing reuse of configurations across differently-sized environments without recomputation. This relates directly to the literature on facility location and emergency asset deployment (\cite{Farahani2009, Murray2013}).

\subsection*{Acknowledgement}
We thank Sahil Wagh, Suvam Das and Sangram, former students of the fourth author at the Department of IE\&OR, IIT Bombay, 
whose early work helped to identify the research direction of this study. We also thank Anjali Bhagat and  Ankita Dargad for many discussions, that helped improve this paper. We thank Abhishek Pathapati, PhD student at the Department of IE\&OR, IIT Bombay, for his help in formulating the integer programming model. Extended thanks to Tirthankar Adhikari whose insights into possible further applications inspired  Section~\ref{sec:thoughts}.  We also gratefully acknowledge the Digital Research Alliance of Canada for providing the computational resources used in this work. Finally, we thank the organizers and participants of Tech Connect 2024 (at IIT Bombay), who ignitioned  this work through a plentitude of interactions between students and attendants of that amazing festival.\\

\end{document}